\documentclass{amsart}
\usepackage{amsmath}
\usepackage{amssymb}
\usepackage{amsfonts}

\newtheorem{theorem}{Theorem}
\theoremstyle{plain}

\newtheorem{corollary}{Corollary}

\newtheorem{definition}{Definition}

\newtheorem{lemma}{Lemma}

\newtheorem{remark}{Remark}

\numberwithin{equation}{section}
\input{tcilatex}

\begin{document}
\title[Some Integral Inequalities]{INEQUALITIES FOR CONVEX AND $s-$CONVEX
FUNCTIONS ON $\Delta =\left[ a,b\right] \times \left[ c,d\right] $}
\author{M. Emin \"{O}zdemir$^{\blacklozenge }$}
\address{$^{\blacklozenge }$Ataturk University, K.K. Education Faculty,
Department of Mathematics, 25240, Erzurum, Turkey}
\email{emos@atauni.edu.tr}
\author{Havva Kavurmac\i $^{\blacklozenge ,\bigstar }$}
\email{hkavurmaci@atauni.edu.tr}
\thanks{$^{\bigstar }$Corresponding Author}
\author{Ahmet Ocak Akdemir$^{\spadesuit }$}
\address{$^{\spadesuit }$A\u{g}r\i\ \.{I}brahim \c{C}e\c{c}en University,
Faculty of Science and Arts, Department of Mathematics, 04100, A\u{g}r\i ,
Turkey}
\email{ahmetakdemir@agri.edu.tr}
\author{Merve Avc\i $^{\blacklozenge }$}
\email{merveavci@ymail.com}
\date{December 29, 2010}

\begin{abstract}
In this paper, two new lemmas are proved and inequalities are established
for co-ordinated convex functions and co-ordinated $s-$convex functions.
\end{abstract}

\maketitle

\section{INTRODUCTION}

Let $f:I\subseteq 
\mathbb{R}
\rightarrow 
\mathbb{R}
$ be a convex function defined on the interval $I$ of real numbers and $a<b.$
The following double inequality;%
\begin{equation*}
f\left( \frac{a+b}{2}\right) \leq \frac{1}{b-a}\dint\limits_{a}^{b}f(x)dx%
\leq \frac{f(a)+f(b)}{2}
\end{equation*}

is well known in the literature as Hadamard's inequality. Both inequalities
hold in the reversed direction if $f$ is concave.

In \cite{OR}, Orlicz defined $s-$convex functions as following:

\begin{definition}
A function $f:%
\mathbb{R}
^{+}\rightarrow 
\mathbb{R}
,$ where $%
\mathbb{R}
^{+}=[0,\infty ),$ is said to be $s-$convex in the first sense if 
\begin{equation*}
f(\alpha x+\beta y)\leq \alpha ^{s}f(x)+\beta ^{s}f(y)
\end{equation*}%
for all $x,y\in \lbrack 0,\infty ),$ $\alpha ,\beta \geq 0$ with $\alpha
^{s}+\beta ^{s}=1$ and for some fixed $s\in (0,1].$ We denote by $K_{s}^{1}$
the class of all $s-$convex functions$.$
\end{definition}

\begin{definition}
A function $f:%
\mathbb{R}
^{+}\rightarrow 
\mathbb{R}
,$ where $%
\mathbb{R}
^{+}=[0,\infty ),$ is said to be $s-$convex in the second sense if 
\begin{equation*}
f(\alpha x+\beta y)\leq \alpha ^{s}f(x)+\beta ^{s}f(y)
\end{equation*}%
for all $x,y\in \lbrack 0,\infty ),$ $\alpha ,\beta \geq 0$ with $\alpha
+\beta =1$ and for some fixed $s\in (0,1].$ We denote by $K_{s}^{2}$ the
class of all $s-$convex functions.
\end{definition}

Obviously, one can see that if we choose $s=1$, both definitions reduced to
ordinary concept of convexity.

For several results related to above definitions we refer readers to \cite%
{HM}, \cite{SF}, \cite{UK}.

In \cite{SS}, Dragomir defined convex functions on the co-ordinates as
following:

\begin{definition}
Let us consider the bidimensional interval $\Delta =[a,b]\times \lbrack c,d]$
in $%
\mathbb{R}
^{2}$ with $a<b,$ $c<d.$ A function $f:\Delta \rightarrow 
\mathbb{R}
$ will be called convex on the co-ordinates if the partial mappings $%
f_{y}:[a,b]\rightarrow 
\mathbb{R}
,$ $f_{y}(u)=f(u,y)$ and $f_{x}:[c,d]\rightarrow 
\mathbb{R}
,$ $f_{x}(v)=f(x,v)$ are convex where defined for all $y\in \lbrack c,d]$
and $x\in \lbrack a,b].$ Recall that the mapping $f:\Delta \rightarrow 
\mathbb{R}
$ is convex on $\Delta $ if the following inequality holds, 
\begin{equation*}
f(\lambda x+(1-\lambda )z,\lambda y+(1-\lambda )w)\leq \lambda
f(x,y)+(1-\lambda )f(z,w)
\end{equation*}%
for all $(x,y),(z,w)\in \Delta $ and $\lambda \in \lbrack 0,1].$
\end{definition}

In \cite{SS}, Dragomir established the following inequalities of Hadamard's
type for co-ordinated convex functions on a rectangle from the plane $%
\mathbb{R}
^{2}.$

\begin{theorem}
Suppose that $f:\Delta =[a,b]\times \lbrack c,d]\rightarrow 
\mathbb{R}
$ is convex on the co-ordinates on $\Delta $. Then one has the inequalities;%
\begin{eqnarray}
&&\ f(\frac{a+b}{2},\frac{c+d}{2})  \label{1.1} \\
&\leq &\frac{1}{2}\left[ \frac{1}{b-a}\int_{a}^{b}f(x,\frac{c+d}{2})dx+\frac{%
1}{d-c}\int_{c}^{d}f(\frac{a+b}{2},y)dy\right]  \notag \\
&\leq &\frac{1}{(b-a)(d-c)}\int_{a}^{b}\int_{c}^{d}f(x,y)dxdy  \notag \\
&\leq &\frac{1}{4}\left[ \frac{1}{(b-a)}\int_{a}^{b}f(x,c)dx+\frac{1}{(b-a)}%
\int_{a}^{b}f(x,d)dx\right.  \notag \\
&&\left. +\frac{1}{(d-c)}\int_{c}^{d}f(a,y)dy+\frac{1}{(d-c)}%
\int_{c}^{d}f(b,y)dy\right]  \notag \\
&\leq &\frac{f(a,c)+f(a,d)+f(b,c)+f(b,d)}{4}.  \notag
\end{eqnarray}%
The above inequalities are sharp.
\end{theorem}

Similar results can be found in \cite{ALO}, \cite{BAK}, \cite{DAR2} and \cite%
{OZ}.

In \cite{DAR2}, Alomari and Darus defined co-ordinated $s-$convex functions
and proved some inequalities based on this definition.

\begin{definition}
Consider the bidimensional interval $\Delta =[a,b]\times \lbrack c,d]$ in $%
[0,\infty )^{2}$ with $a<b$ and $c<d$. The mapping $f:\Delta \rightarrow 
\mathbb{R}
$ is $s-$convex on $\Delta $ if 
\begin{equation*}
f(\lambda x+(1-\lambda )z,\lambda y+(1-\lambda )w)\leq \lambda
^{s}f(x,y)+(1-\lambda )^{s}f(z,w)
\end{equation*}%
holds for all $(x,y),$ $(z,w)\in \Delta $ with $\lambda \in \lbrack 0,1]$
and for some fixed $s\in (0,1].$
\end{definition}

In \cite{OZ2}, Sar\i kaya \textit{et al.} proved some Hadamard's type
inequalities for co-ordinated convex functions as followings:

\begin{theorem}
Let $f:\Delta \subset 
\mathbb{R}
^{2}\rightarrow 
\mathbb{R}
$ be a partial differentiable mapping on $\Delta :=[a,b]\times \lbrack c,d]$
in $%
\mathbb{R}
^{2}$ with $a<b$ and $c<d.$ If $\left\vert \frac{\partial ^{2}f}{\partial
t\partial s}\right\vert $ is a convex function on the co-ordinates on $%
\Delta ,$ then one has the inequalities:%
\begin{eqnarray}
&&\left\vert \frac{f(a,c)+f(a,d)+f(b,c)+f(b,d)}{4}\right.  \label{1.3} \\
&&\left. \frac{1}{(b-a)(d-c)}\int_{a}^{b}\int_{c}^{d}f(x,y)dxdy-A\right\vert
\notag \\
&\leq &\frac{(b-a)(d-c)}{16}\left( \frac{\left\vert \frac{\partial ^{2}f}{%
\partial t\partial s}\right\vert (a,c)+\left\vert \frac{\partial ^{2}f}{%
\partial t\partial s}\right\vert (a,d)+\left\vert \frac{\partial ^{2}f}{%
\partial t\partial s}\right\vert (b,c)+\left\vert \frac{\partial ^{2}f}{%
\partial t\partial s}\right\vert (b,d)}{4}\right)  \notag
\end{eqnarray}%
where%
\begin{equation*}
A=\frac{1}{2}\left[ \frac{1}{(b-a)}\int_{a}^{b}\left[ f(x,c)+f(x,d)\right]
dx+\frac{1}{(d-c)}\int_{c}^{d}\left[ f(a,y)dy+f(b,y)\right] dy\right] .
\end{equation*}
\end{theorem}

\begin{theorem}
Let $f:\Delta \subset 
\mathbb{R}
^{2}\rightarrow 
\mathbb{R}
$ be a partial differentiable mapping on $\Delta :=[a,b]\times \lbrack c,d]$
in $%
\mathbb{R}
^{2}$ with $a<b$ and $c<d.$ If $\left\vert \frac{\partial ^{2}f}{\partial
t\partial s}\right\vert ^{q},$ $q>1,$ is a convex function on the
co-ordinates on $\Delta ,$ then one has the inequalities:%
\begin{eqnarray}
&&\left\vert \frac{f(a,c)+f(a,d)+f(b,c)+f(b,d)}{4}\right.  \label{1.4} \\
&&\left. \frac{1}{(b-a)(d-c)}\int_{a}^{b}\int_{c}^{d}f(x,y)dxdy-A\right\vert
\notag \\
&\leq &\frac{(b-a)(d-c)}{4\left( p+1\right) ^{\frac{2}{p}}}\left( \frac{%
\left\vert \frac{\partial ^{2}f}{\partial t\partial s}\right\vert
^{q}(a,c)+\left\vert \frac{\partial ^{2}f}{\partial t\partial s}\right\vert
^{q}(a,d)+\left\vert \frac{\partial ^{2}f}{\partial t\partial s}\right\vert
^{q}(b,c)+\left\vert \frac{\partial ^{2}f}{\partial t\partial s}\right\vert
^{q}(b,d)}{4}\right) ^{\frac{1}{q}}  \notag
\end{eqnarray}%
where%
\begin{equation*}
A=\frac{1}{2}\left[ \frac{1}{(b-a)}\int_{a}^{b}\left[ f(x,c)+f(x,d)\right]
dx+\frac{1}{(d-c)}\int_{c}^{d}\left[ f(a,y)dy+f(b,y)\right] dy\right]
\end{equation*}%
and $\frac{1}{p}+\frac{1}{q}=1.$
\end{theorem}

\begin{theorem}
Let $f:\Delta \subset 
\mathbb{R}
^{2}\rightarrow 
\mathbb{R}
$ be a partial differentiable mapping on $\Delta :=[a,b]\times \lbrack c,d]$
in $%
\mathbb{R}
^{2}$ with $a<b$ and $c<d.$ If $\left\vert \frac{\partial ^{2}f}{\partial
t\partial s}\right\vert ^{q},$ $q\geq 1,$ is a convex function on the
co-ordinates on $\Delta ,$ then one has the inequalities:%
\begin{eqnarray}
&&\left\vert \frac{f(a,c)+f(a,d)+f(b,c)+f(b,d)}{4}\right.  \label{1.5} \\
&&\left. \frac{1}{(b-a)(d-c)}\int_{a}^{b}\int_{c}^{d}f(x,y)dxdy-A\right\vert
\notag \\
&\leq &\frac{(b-a)(d-c)}{16}\left( \frac{\left\vert \frac{\partial ^{2}f}{%
\partial t\partial s}\right\vert ^{q}(a,c)+\left\vert \frac{\partial ^{2}f}{%
\partial t\partial s}\right\vert ^{q}(a,d)+\left\vert \frac{\partial ^{2}f}{%
\partial t\partial s}\right\vert ^{q}(b,c)+\left\vert \frac{\partial ^{2}f}{%
\partial t\partial s}\right\vert ^{q}(b,d)}{4}\right) ^{\frac{1}{q}}  \notag
\end{eqnarray}%
where%
\begin{equation*}
A=\frac{1}{2}\left[ \frac{1}{(b-a)}\int_{a}^{b}\left[ f(x,c)+f(x,d)\right]
dx+\frac{1}{(d-c)}\int_{c}^{d}\left[ f(a,y)dy+f(b,y)\right] dy\right] .
\end{equation*}
\end{theorem}

In \cite{SS2}, Barnett and Dragomir proved an Ostrowski-type inequality for
double integrals as following:

\begin{theorem}
\label{1} Let $f:\left[ a,b\right] \times \left[ c,d\right] \rightarrow 
\mathbb{R}
$ be continuous on $\left[ a,b\right] \times \left[ c,d\right] $, $%
f_{x,y}^{\prime \prime }=\frac{\partial ^{2}f}{\partial x\partial y}$ exists
on $\left( a,b\right) \times \left( c,d\right) $ and is bounded, that is%
\begin{equation*}
\left\Vert f_{x,y}^{\prime \prime }\right\Vert _{\infty }=\sup_{\left(
x,y\right) \in \left( a,b\right) \times \left( c,d\right) }\left\vert \frac{%
\partial ^{2}f\left( x,y\right) }{\partial x\partial y}\right\vert <\infty ,
\end{equation*}%
then we have the inequality;%
\begin{eqnarray}
&&\left\vert
\dint\limits_{a}^{b}\dint\limits_{c}^{d}f(s,t)dtds-(b-a)\dint%
\limits_{c}^{d}f(x,t)dt\right.  \label{1.6} \\
&&\left. -(d-c)\dint\limits_{a}^{b}f(s,y)ds-(b-a)(d-c)f(x,y)\right\vert 
\notag \\
&\leq &\left[ \frac{(b-a)^{2}}{4}+\left( x-\frac{a+b}{2}\right) ^{2}\right] %
\left[ \frac{(d-c)^{2}}{4}+\left( y-\frac{c+d}{2}\right) ^{2}\right]
\left\Vert f_{x,y}^{\prime \prime }\right\Vert _{\infty }  \notag
\end{eqnarray}%
for all $(x,y)\in \left[ a,b\right] \times \left[ c,d\right] $.
\end{theorem}

In \cite{ZEKI2}, Sar\i kaya proved an Ostrowski-type inequality for double
integrals and gave a corollary as following:

\begin{theorem}
\label{2} Let $f:\left[ a,b\right] \times \left[ c,d\right] \rightarrow 
\mathbb{R}
$ be an absolutely continuous functions such that the partial derivative of
order 2 exist and is bounded, i.e.,%
\begin{equation*}
\left\Vert \frac{\partial ^{2}f\left( t,s\right) }{\partial t\partial s}%
\right\Vert _{\infty }=\underset{\left( x,y\right) \in \left( a,b\right)
\times \left( c,d\right) }{\sup }\left\vert \frac{\partial ^{2}f\left(
t,s\right) }{\partial t\partial s}\right\vert <\infty
\end{equation*}%
for all $\left( t,s\right) \in \left[ a,b\right] \times \left[ c,d\right] $.
Then we have,%
\begin{eqnarray}
&&\left\vert \left( \beta _{1}-\alpha _{1}\right) \left( \beta _{2}-\alpha
_{2}\right) f\left( \frac{a+b}{2},\frac{c+d}{2}\right) +H\left( \alpha
_{1},\alpha _{2},\beta _{1},\beta _{2}\right) +G\left( \alpha _{1},\alpha
_{2},\beta _{1},\beta _{2}\right) \right.  \notag \\
&&-\left( \beta _{2}-\alpha _{2}\right) \int_{a}^{b}f\left( t,\frac{c+d}{2}%
\right) dt-\left( \beta _{1}-\alpha _{1}\right) \int_{c}^{d}f\left( \frac{a+b%
}{2},s\right) ds  \notag \\
&&-\int_{a}^{b}\left[ \left( \alpha _{2}-c\right) f\left( t,c\right) +\left(
d-\beta _{2}\right) f(t,d)\right] dt  \label{1.7} \\
&&\left. -\int_{c}^{d}\left[ \left( \alpha _{1}-a\right) f\left( a,s\right)
+\left( b-\beta _{1}\right) f(b,s)\right] ds+\int_{a}^{b}%
\int_{c}^{d}f(t,s)dsdt\right\vert  \notag \\
&\leq &\left[ \frac{\left( \alpha _{1}-a\right) ^{2}+\left( b-\beta
_{1}\right) ^{2}}{2}+\frac{\left( a+b-2\alpha _{1}\right) ^{2}+\left(
a+b-2\beta _{1}\right) ^{2}}{8}\right]  \notag \\
&&\times \left[ \frac{\left( \alpha _{2}-c\right) ^{2}+\left( d-\beta
_{2}\right) ^{2}}{2}+\frac{\left( c+d-2\alpha _{2}\right) ^{2}+\left(
c+d-2\beta _{2}\right) ^{2}}{8}\right] \left\Vert \frac{\partial ^{2}f\left(
t,s\right) }{\partial t\partial s}\right\Vert _{\infty }  \notag
\end{eqnarray}%
for all $\left( \alpha _{1},\alpha _{2}\right) ,\left( \beta _{1},\beta
_{2}\right) \in $ $\left[ a,b\right] \times \left[ c,d\right] $ with $\alpha
_{1}<\beta _{1},$ $\alpha _{2}<\beta _{2}$ where%
\begin{eqnarray*}
&&H\left( \alpha _{1},\alpha _{2},\beta _{1},\beta _{2}\right) \\
&=&\left( \alpha _{1}-a\right) \left[ \left( \alpha _{2}-c\right)
f(a,c)+\left( d-\beta _{2}\right) f(a,d)\right] \\
&&+\left( b-\beta _{1}\right) \left[ \left( \alpha _{2}-c\right)
f(b,c)+(d-\beta _{2})f(b,d)\right]
\end{eqnarray*}%
and%
\begin{eqnarray*}
&&G\left( \alpha _{1},\alpha _{2},\beta _{1},\beta _{2}\right) \\
&=&\left( \beta _{1}-\alpha _{1}\right) \left[ \left( \alpha _{2}-c\right)
f\left( \frac{a+b}{2},c\right) +\left( d-\beta _{2}\right) f\left( \frac{a+b%
}{2},d\right) \right] \\
&&+\left( \beta _{2}-\alpha _{2}\right) \left[ \left( \alpha _{1}-a\right)
f\left( a,\frac{c+d}{2}\right) +\left( b-\beta _{1}\right) f\left( b,\frac{%
c+d}{2}\right) \right] .
\end{eqnarray*}
\end{theorem}

\begin{corollary}
Under the assumptions of Theorem \ref{2}, we have%
\begin{eqnarray}
&&\left\vert \left( b-a\right) \left( d-c\right) f\left( \frac{a+b}{2},\frac{%
c+d}{2}\right) +\int_{a}^{b}\int_{c}^{d}f(t,s)dsdt\right.  \label{1.8} \\
&&\left. -\left( d-c\right) \int_{a}^{b}f\left( t,\frac{c+d}{2}\right)
dt-\left( b-a\right) \int_{c}^{d}f\left( \frac{a+b}{2},s\right) ds\right\vert
\notag \\
&\leq &\frac{1}{16}\left\Vert \frac{\partial ^{2}f\left( t,s\right) }{%
\partial t\partial s}\right\Vert _{\infty }\left( b-a\right) ^{2}\left(
d-c\right) ^{2}.  \notag
\end{eqnarray}
\end{corollary}

In \cite{PA}, Pachpatte established a new Ostrowski type inequality similar
to inequality (\ref{1.6}) by using elemantery analysis.

The main purpose of this paper is to establish inequalities of
Ostrowski-type for co-ordinated convex functions by using Lemma 2 and to
establish some new Hadamard's type inequalities for co-ordinated $s-$convex
functions.

\section{INEQUALITIES\ FOR CO-ORDINATED CONVEX FUNCTIONS}

To prove our main result, we need the following lemma which contains kernels
similar to Barnett and Dragomir's kernels in \cite{SS2}, [see the paper \cite%
{SS2}, proof of Theorem 2.1].

\begin{lemma}
Let $f:\Delta =\left[ a,b\right] \times \left[ c,d\right] \rightarrow 
\mathbb{R}
$ be a partial differentiable mapping on $\Delta =\left[ a,b\right] \times %
\left[ c,d\right] .$ If $\frac{\partial ^{2}f}{\partial t\partial s}\in
L\left( \Delta \right) ,$ then the following equality holds:%
\begin{eqnarray*}
&&f\left( \frac{a+b}{2},\frac{c+d}{2}\right) \\
&&-\frac{1}{\left( d-c\right) }\int_{c}^{d}f\left( \frac{a+b}{2},y\right) dy-%
\frac{1}{\left( b-a\right) }\int_{a}^{b}f\left( x,\frac{c+d}{2}\right) dx \\
&&+\frac{1}{\left( b-a\right) \left( d-c\right) }\int_{a}^{b}\int_{c}^{d}f%
\left( x,y\right) dydx \\
&=&\frac{1}{\left( b-a\right) \left( d-c\right) }\int_{a}^{b}\int_{c}^{d}p%
\left( x,t\right) q\left( y,s\right) \frac{\partial ^{2}f}{\partial
t\partial s}\left( \frac{b-t}{b-a}a+\frac{t-a}{b-a}b,\frac{d-s}{d-c}c+\frac{%
s-c}{d-c}d\right) dsdt
\end{eqnarray*}%
where%
\begin{equation*}
p\left( x,t\right) =\left\{ 
\begin{array}{c}
\left( t-a\right) ,\text{ \ \ \ \ }t\in \left[ a,\frac{a+b}{2}\right] \\ 
\\ 
\left( t-b\right) ,\text{ \ \ \ \ }t\in \left( \frac{a+b}{2},b\right]%
\end{array}%
\right.
\end{equation*}%
and%
\begin{equation*}
q\left( y,s\right) =\left\{ 
\begin{array}{c}
\left( s-c\right) ,\text{ \ \ \ \ }s\in \left[ c,\frac{c+d}{2}\right] \\ 
\\ 
\left( s-d\right) ,\text{ \ \ \ \ }s\in \left( \frac{c+d}{2},d\right]%
\end{array}%
\right. .
\end{equation*}%
for each $x\in \left[ a,b\right] $ and $y\in \left[ c,d\right] .$
\end{lemma}

\begin{proof}
Integration by parts, we can write%
\begin{eqnarray*}
&&\frac{1}{\left( b-a\right) \left( d-c\right) }\int_{a}^{b}\int_{c}^{d}p%
\left( x,t\right) q\left( y,s\right) \frac{\partial ^{2}f}{\partial
t\partial s}\left( \frac{b-t}{b-a}a+\frac{t-a}{b-a}b,\frac{d-s}{d-c}c+\frac{%
s-c}{d-c}d\right) dsdt \\
&=&\int_{c}^{d}q\left( y,s\right) \left[ \int_{a}^{\frac{a+b}{2}}\left(
t-a\right) \frac{\partial ^{2}f}{\partial t\partial s}\left( \frac{b-t}{b-a}%
a+\frac{t-a}{b-a}b,\frac{d-s}{d-c}c+\frac{s-c}{d-c}d\right) dt\right. \\
&&\left. +\int_{\frac{a+b}{2}}^{b}\left( t-b\right) \frac{\partial ^{2}f}{%
\partial t\partial s}\left( \frac{b-t}{b-a}a+\frac{t-a}{b-a}b,\frac{d-s}{d-c}%
c+\frac{s-c}{d-c}d\right) dt\right] ds
\end{eqnarray*}%
\begin{eqnarray*}
&=&\int_{c}^{d}q\left( y,s\right) \left\{ \left[ \left( t-a\right) \frac{%
\partial f}{\partial s}\left( \frac{b-t}{b-a}a+\frac{t-a}{b-a}b,\frac{d-s}{%
d-c}c+\frac{s-c}{d-c}d\right) \right] _{a}^{\frac{a+b}{2}}\right. \\
&&-\int_{a}^{\frac{a+b}{2}}\frac{\partial f}{\partial s}\left( \frac{b-t}{b-a%
}a+\frac{t-a}{b-a}b,\frac{d-s}{d-c}c+\frac{s-c}{d-c}d\right) dt \\
&&+\left[ \left( t-b\right) \frac{\partial f}{\partial s}\left( \frac{b-t}{%
b-a}a+\frac{t-a}{b-a}b,\frac{d-s}{d-c}c+\frac{s-c}{d-c}d\right) \right] _{%
\frac{a+b}{2}}^{b} \\
&&\left. -\int_{\frac{a+b}{2}}^{b}\frac{\partial f}{\partial s}\left( \frac{%
b-t}{b-a}a+\frac{t-a}{b-a}b,\frac{d-s}{d-c}c+\frac{s-c}{d-c}d\right)
dt\right\} ds.
\end{eqnarray*}%
We obtain%
\begin{eqnarray*}
&&\frac{1}{\left( b-a\right) \left( d-c\right) }\int_{a}^{b}\int_{c}^{d}p%
\left( x,t\right) q\left( y,s\right) \frac{\partial ^{2}f}{\partial
t\partial s}\left( \frac{b-t}{b-a}a+\frac{t-a}{b-a}b,\frac{d-s}{d-c}c+\frac{%
s-c}{d-c}d\right) dsdt \\
&=&\left( b-a\right) \int_{c}^{d}q\left( y,s\right) \left\{ \frac{\partial f%
}{\partial s}\left( \frac{a+b}{2},\frac{d-s}{d-c}c+\frac{s-c}{d-c}d\right)
\right. \\
&&\left. -\int_{a}^{b}\frac{\partial f}{\partial s}\left( \frac{b-t}{b-a}a+%
\frac{t-a}{b-a}b,\frac{d-s}{d-c}c+\frac{s-c}{d-c}d\right) dt\right\} ds.
\end{eqnarray*}%
By integrating again, we get%
\begin{eqnarray*}
&&\frac{1}{\left( b-a\right) \left( d-c\right) }\int_{a}^{b}\int_{c}^{d}p%
\left( x,t\right) q\left( y,s\right) \frac{\partial ^{2}f}{\partial
t\partial s}\left( \frac{b-t}{b-a}a+\frac{t-a}{b-a}b,\frac{d-s}{d-c}c+\frac{%
s-c}{d-c}d\right) dsdt \\
&=&\left( b-a\right) \left\{ \int_{c}^{\frac{c+d}{2}}\left( s-c\right) \frac{%
\partial f}{\partial s}\left( \frac{a+b}{2},\frac{d-s}{d-c}c+\frac{s-c}{d-c}%
d\right) ds\right. \\
&&+\int_{\frac{c+d}{2}}^{d}\left( s-d\right) \frac{\partial f}{\partial s}%
\left( \frac{a+b}{2},\frac{d-s}{d-c}c+\frac{s-c}{d-c}d\right) ds \\
&&-\int_{a}^{b}\left[ \int_{c}^{\frac{c+d}{2}}\left( s-c\right) \frac{%
\partial f}{\partial s}\left( \frac{b-t}{b-a}a+\frac{t-a}{b-a}b,\frac{d-s}{%
d-c}c+\frac{s-c}{d-c}d\right) ds\right. \\
&&\left. \left. +\int_{\frac{c+d}{2}}^{d}\left( s-d\right) \frac{\partial f}{%
\partial s}\left( \frac{b-t}{b-a}a+\frac{t-a}{b-a}b,\frac{d-s}{d-c}c+\frac{%
s-c}{d-c}d\right) ds\right] dt\right\} .
\end{eqnarray*}%
By calculating the above integrals, we have%
\begin{eqnarray*}
&&\frac{1}{\left( b-a\right) \left( d-c\right) }\int_{a}^{b}\int_{c}^{d}p%
\left( x,t\right) q\left( y,s\right) \frac{\partial ^{2}f}{\partial
t\partial s}\left( \frac{b-t}{b-a}a+\frac{t-a}{b-a}b,\frac{d-s}{d-c}c+\frac{%
s-c}{d-c}d\right) dsdt \\
&=&\left( b-a\right) \left( d-c\right) f\left( \frac{a+b}{2},\frac{c+d}{2}%
\right) \\
&&-\left( b-a\right) \int_{c}^{d}f\left( \frac{a+b}{2},\frac{d-s}{d-c}c+%
\frac{s-c}{d-c}d\right) ds \\
&&-\left( d-c\right) \int_{a}^{b}f\left( \frac{b-t}{b-a}a+\frac{t-a}{b-a}b,%
\frac{c+d}{2}\right) dt \\
&&\int_{a}^{b}\int_{c}^{d}f\left( \frac{b-t}{b-a}a+\frac{t-a}{b-a}b,\frac{d-s%
}{d-c}c+\frac{s-c}{d-c}d\right) dsdt.
\end{eqnarray*}%
Using the change of the variable $x=\frac{b-t}{b-a}a+\frac{t-a}{b-a}b$ and $%
y=\frac{d-s}{d-c}c+\frac{s-c}{d-c}d,$ then dividing both sides with $\left(
b-a\right) \times \left( d-c\right) ,$ this completes the proof.
\end{proof}

\begin{theorem}
\label{t.2.1} Let $f:\Delta =\left[ a,b\right] \times \left[ c,d\right]
\rightarrow 
\mathbb{R}
$ be a partial differentiable mapping on $\Delta =\left[ a,b\right] \times %
\left[ c,d\right] .$ If $\left\vert \frac{\partial ^{2}f}{\partial t\partial
s}\right\vert $ is a convex function on the co-ordinates on $\Delta ,$ then
the following inequality holds;%
\begin{eqnarray}
&&\left\vert f\left( \frac{a+b}{2},\frac{c+d}{2}\right) \right.  \label{2.1}
\\
&&-\frac{1}{\left( d-c\right) }\int_{c}^{d}f\left( \frac{a+b}{2},y\right) dy-%
\frac{1}{\left( b-a\right) }\int_{a}^{b}f\left( x,\frac{c+d}{2}\right) dx 
\notag \\
&&\left. +\frac{1}{\left( b-a\right) \left( d-c\right) }\int_{a}^{b}%
\int_{c}^{d}f\left( x,y\right) dydx\right\vert  \notag \\
&\leq &\frac{\left( b-a\right) \left( d-c\right) }{64}\left[ \left\vert 
\frac{\partial ^{2}f}{\partial t\partial s}\left( a,c\right) \right\vert
+\left\vert \frac{\partial ^{2}f}{\partial t\partial s}\left( b,c\right)
\right\vert +\left\vert \frac{\partial ^{2}f}{\partial t\partial s}\left(
a,d\right) \right\vert +\left\vert \frac{\partial ^{2}f}{\partial t\partial s%
}\left( b,d\right) \right\vert \right] .  \notag
\end{eqnarray}
\end{theorem}

\begin{proof}
From Lemma 2 and using the property of modulus, we have%
\begin{eqnarray*}
&&\left\vert f\left( \frac{a+b}{2},\frac{c+d}{2}\right) \right. \\
&&-\frac{1}{\left( d-c\right) }\int_{c}^{d}f\left( \frac{a+b}{2},y\right) dy-%
\frac{1}{\left( b-a\right) }\int_{a}^{b}f\left( x,\frac{c+d}{2}\right) dx \\
&&\left. +\frac{1}{\left( b-a\right) \left( d-c\right) }\int_{a}^{b}%
\int_{c}^{d}f\left( x,y\right) dydx\right\vert \\
&\leq &\frac{1}{\left( b-a\right) \left( d-c\right) }\int_{a}^{b}%
\int_{c}^{d}\left\vert p\left( x,t\right) q\left( y,s\right) \right\vert
\left\vert \frac{\partial ^{2}f}{\partial t\partial s}\left( \frac{b-t}{b-a}%
a+\frac{t-a}{b-a}b,\frac{d-s}{d-c}c+\frac{s-c}{d-c}d\right) \right\vert dsdt
\end{eqnarray*}%
Since $\left\vert \frac{\partial ^{2}f}{\partial t\partial s}\right\vert $
is co-ordinated convex, we can write%
\begin{eqnarray*}
&&\left\vert f\left( \frac{a+b}{2},\frac{c+d}{2}\right) \right. \\
&&-\frac{1}{\left( d-c\right) }\int_{c}^{d}f\left( \frac{a+b}{2},y\right) dy-%
\frac{1}{\left( b-a\right) }\int_{a}^{b}f\left( x,\frac{c+d}{2}\right) dx \\
&&\left. +\frac{1}{\left( b-a\right) \left( d-c\right) }\int_{a}^{b}%
\int_{c}^{d}f\left( x,y\right) dydx\right\vert \\
&\leq &\frac{1}{\left( b-a\right) \left( d-c\right) } \\
&&\times \int_{c}^{d}\left\vert q\left( y,s\right) \right\vert \left\{
\int_{a}^{\frac{a+b}{2}}\left( t-a\right) \left[ \frac{b-t}{b-a}\left\vert 
\frac{\partial ^{2}f}{\partial t\partial s}\left( a,\frac{d-s}{d-c}c+\frac{%
s-c}{d-c}d\right) \right\vert \right] dt\right. \\
&&+\int_{a}^{\frac{a+b}{2}}\left( t-a\right) \left[ \frac{t-a}{b-a}%
\left\vert \frac{\partial ^{2}f}{\partial t\partial s}\left( b,\frac{d-s}{d-c%
}c+\frac{s-c}{d-c}d\right) \right\vert \right] dt \\
&&+\int_{\frac{a+b}{2}}^{b}\left( b-t\right) \left[ \frac{b-t}{b-a}%
\left\vert \frac{\partial ^{2}f}{\partial t\partial s}\left( a,\frac{d-s}{d-c%
}c+\frac{s-c}{d-c}d\right) \right\vert \right] dt \\
&&\left. +\int_{\frac{a+b}{2}}^{b}\left( b-t\right) \left[ \frac{t-a}{b-a}%
\left\vert \frac{\partial ^{2}f}{\partial t\partial s}\left( b,\frac{d-s}{d-c%
}c+\frac{s-c}{d-c}d\right) \right\vert \right] dt\right\} ds.
\end{eqnarray*}%
By computing these integrals, we obtain%
\begin{eqnarray*}
&&\left\vert f\left( \frac{a+b}{2},\frac{c+d}{2}\right) \right. \\
&&-\frac{1}{\left( d-c\right) }\int_{c}^{d}f\left( \frac{a+b}{2},y\right) ds-%
\frac{1}{\left( b-a\right) }\int_{a}^{b}f\left( x,\frac{c+d}{2}\right) dx \\
&&\left. +\frac{1}{\left( b-a\right) \left( d-c\right) }\int_{a}^{b}%
\int_{c}^{d}f\left( x,y\right) dxdy\right\vert \\
&\leq &\frac{\left( b-a\right) }{8\left( d-c\right) }\left[
\int_{c}^{d}\left\vert q\left( y,s\right) \right\vert \left\vert \frac{%
\partial ^{2}f}{\partial t\partial s}\left( a,\frac{d-s}{d-c}c+\frac{s-c}{d-c%
}d\right) \right\vert \right. \\
&&+\left. \int_{c}^{d}\left\vert q\left( y,s\right) \right\vert \left\vert 
\frac{\partial ^{2}f}{\partial t\partial s}\left( b,\frac{d-s}{d-c}c+\frac{%
s-c}{d-c}d\right) \right\vert \right] ds.
\end{eqnarray*}%
Using co-ordinated convexity of $\left\vert \frac{\partial ^{2}f}{\partial
t\partial s}\right\vert $ again, we get%
\begin{eqnarray*}
&&\left\vert f\left( \frac{a+b}{2},\frac{c+d}{2}\right) \right. \\
&&-\frac{1}{\left( d-c\right) }\int_{c}^{d}f\left( \frac{a+b}{2},y\right) dy-%
\frac{1}{\left( b-a\right) }\int_{a}^{b}f\left( x,\frac{c+d}{2}\right) dx \\
&&\left. +\frac{1}{\left( b-a\right) \left( d-c\right) }\int_{a}^{b}%
\int_{c}^{d}f\left( x,y\right) dydx\right\vert \\
&\leq &\frac{\left( b-a\right) }{8\left( d-c\right) }\left[ \int_{c}^{\frac{%
c+d}{2}}\left( s-c\right) \left[ \frac{d-s}{d-c}\left\vert \frac{\partial
^{2}f}{\partial t\partial s}\left( a,c\right) \right\vert \right]
ds+\int_{c}^{\frac{c+d}{2}}\left( s-c\right) \left[ \frac{s-c}{d-c}%
\left\vert \frac{\partial ^{2}f}{\partial t\partial s}\left( a,d\right)
\right\vert \right] ds\right. \\
&&+\int_{\frac{c+d}{2}}^{d}\left( d-s\right) \left[ \frac{d-s}{d-c}%
\left\vert \frac{\partial ^{2}f}{\partial t\partial s}\left( a,c\right)
\right\vert \right] ds+\int_{\frac{c+d}{2}}^{d}\left( d-s\right) \left[ 
\frac{s-c}{d-c}\left\vert \frac{\partial ^{2}f}{\partial t\partial s}\left(
a,d\right) \right\vert \right] ds \\
&&+\int_{c}^{\frac{c+d}{2}}\left( s-c\right) \left[ \frac{d-s}{d-c}%
\left\vert \frac{\partial ^{2}f}{\partial t\partial s}\left( b,c\right)
\right\vert \right] ds+\int_{c}^{\frac{c+d}{2}}\left( s-c\right) \left[ 
\frac{s-c}{d-c}\left\vert \frac{\partial ^{2}f}{\partial t\partial s}\left(
b,d\right) \right\vert \right] ds \\
&&\left. +\int_{\frac{c+d}{2}}^{d}\left( d-s\right) \left[ \frac{d-s}{d-c}%
\left\vert \frac{\partial ^{2}f}{\partial t\partial s}\left( b,c\right)
\right\vert \right] ds+\int_{\frac{c+d}{2}}^{d}\left( d-s\right) \left[ 
\frac{s-c}{d-c}\left\vert \frac{\partial ^{2}f}{\partial t\partial s}\left(
b,d\right) \right\vert \right] ds\right\} .
\end{eqnarray*}%
By a simple computation, we get the required result.
\end{proof}

\begin{remark}
Suppose that all the assumptions of Theorem \ref{t.2.1} are satisfied. If we
choose $\frac{\partial ^{2}f}{\partial t\partial s}$ is bounded, i.e.,%
\begin{equation*}
\left\Vert \frac{\partial ^{2}f\left( t,s\right) }{\partial t\partial s}%
\right\Vert _{\infty }=\underset{\left( t,s\right) \in \left( a,b\right)
\times \left( c,d\right) }{\sup }\left\vert \frac{\partial ^{2}f\left(
t,s\right) }{\partial t\partial s}\right\vert <\infty ,
\end{equation*}%
we get%
\begin{eqnarray}
&&\left\vert f\left( \frac{a+b}{2},\frac{c+d}{2}\right) \right.
\label{r.2.1} \\
&&-\frac{1}{\left( d-c\right) }\int_{c}^{d}f\left( \frac{a+b}{2},y\right) dy-%
\frac{1}{\left( b-a\right) }\int_{a}^{b}f\left( x,\frac{c+d}{2}\right) dx 
\notag \\
&&\left. +\frac{1}{\left( b-a\right) \left( d-c\right) }\int_{a}^{b}%
\int_{c}^{d}f\left( x,y\right) dydx\right\vert  \notag \\
&\leq &\frac{\left( b-a\right) \left( d-c\right) }{16}\left\Vert \frac{%
\partial ^{2}f\left( t,s\right) }{\partial t\partial s}\right\Vert _{\infty }
\notag
\end{eqnarray}%
which is the inequality $\left( \ref{1.8}\right) .$
\end{remark}

\begin{theorem}
\label{t.2.2} Let $f:\Delta =\left[ a,b\right] \times \left[ c,d\right]
\rightarrow 
\mathbb{R}
$ be a partial differentiable mapping on $\Delta =\left[ a,b\right] \times %
\left[ c,d\right] .$ If $\left\vert \frac{\partial ^{2}f}{\partial t\partial
s}\right\vert ^{q},$ $q>1,$ is a convex function on the co-ordinates on $%
\Delta ,$ then the following inequality holds;%
\begin{eqnarray}
&&\left\vert f\left( \frac{a+b}{2},\frac{c+d}{2}\right) +\frac{1}{\left(
b-a\right) \left( d-c\right) }\int_{a}^{b}\int_{c}^{d}f\left( x,y\right)
dydx\right.  \label{2.2} \\
&&\left. -\frac{1}{\left( d-c\right) }\int_{c}^{d}f\left( \frac{a+b}{2}%
,y\right) dy-\frac{1}{\left( b-a\right) }\int_{a}^{b}f\left( x,\frac{c+d}{2}%
\right) dx\right\vert  \notag \\
&\leq &\frac{\left( b-a\right) \left( d-c\right) }{4\left( p+1\right) ^{%
\frac{2}{p}}}  \notag \\
&&\times \left( \frac{\left\vert \frac{\partial ^{2}f}{\partial t\partial s}%
\left( a,c\right) \right\vert ^{q}+\left\vert \frac{\partial ^{2}f}{\partial
t\partial s}\left( b,c\right) \right\vert ^{q}+\left\vert \frac{\partial
^{2}f}{\partial t\partial s}\left( a,d\right) \right\vert ^{q}+\left\vert 
\frac{\partial ^{2}f}{\partial t\partial s}\left( b,d\right) \right\vert ^{q}%
}{4}\right) ^{\frac{1}{q}}.  \notag
\end{eqnarray}
\end{theorem}

\begin{proof}
From Lemma 2, we have%
\begin{eqnarray*}
&&\left\vert f\left( \frac{a+b}{2},\frac{c+d}{2}\right) +\frac{1}{\left(
b-a\right) \left( d-c\right) }\int_{a}^{b}\int_{c}^{d}f\left( x,y\right)
dydx\right. \\
&&\left. -\frac{1}{\left( d-c\right) }\int_{c}^{d}f\left( \frac{a+b}{2}%
,y\right) dy-\frac{1}{\left( b-a\right) }\int_{a}^{b}f\left( x,\frac{c+d}{2}%
\right) dx\right\vert \\
&\leq &\frac{1}{\left( b-a\right) \left( d-c\right) }\int_{a}^{b}%
\int_{c}^{d}\left\vert p\left( x,t\right) q\left( y,s\right) \right\vert
\left\vert \frac{\partial ^{2}f}{\partial t\partial s}\left( \frac{b-t}{b-a}%
a+\frac{t-a}{b-a}b,\frac{d-s}{d-c}c+\frac{s-c}{d-c}d\right) \right\vert dsdt
\end{eqnarray*}%
By applying the well-known H\"{o}lder inequality for double integrals, then
one has%
\begin{eqnarray}
&&\left\vert f\left( \frac{a+b}{2},\frac{c+d}{2}\right) +\frac{1}{\left(
b-a\right) \left( d-c\right) }\int_{a}^{b}\int_{c}^{d}f\left( x,y\right)
dydx\right.  \label{2.3} \\
&&\left. -\frac{1}{\left( d-c\right) }\int_{c}^{d}f\left( \frac{a+b}{2}%
,y\right) dy-\frac{1}{\left( b-a\right) }\int_{a}^{b}f\left( x,\frac{c+d}{2}%
\right) dx\right\vert  \notag \\
&\leq &\frac{1}{\left( b-a\right) \left( d-c\right) }\left\{ \left(
\int_{a}^{b}\int_{c}^{d}\left\vert p\left( x,t\right) q\left( y,s\right)
\right\vert ^{p}dtds\right) ^{\frac{1}{p}}\right.  \notag \\
&&\times \left( \int_{a}^{b}\int_{c}^{d}\left\vert \frac{\partial ^{2}f}{%
\partial t\partial s}\left( \frac{b-t}{b-a}a+\frac{t-a}{b-a}b,\frac{d-s}{d-c}%
c+\frac{s-c}{d-c}d\right) \right\vert ^{q}dsdt\right) ^{\frac{1}{q}}.  \notag
\end{eqnarray}%
Since $\left\vert \frac{\partial ^{2}f}{\partial t\partial s}\right\vert
^{q} $ is a co-ordinated convex function on $\Delta ,$ we can write for all $%
\left( t,s\right) \in \left[ a,b\right] \times \left[ c,d\right] $%
\begin{eqnarray*}
&&\left\vert \frac{\partial ^{2}f}{\partial t\partial s}\left( \frac{b-t}{b-a%
}a+\frac{t-a}{b-a}b,\frac{d-s}{d-c}c+\frac{s-c}{d-c}d\right) \right\vert ^{q}
\\
&\leq &\frac{b-t}{b-a}\left\vert \frac{\partial ^{2}f}{\partial t\partial s}%
\left( a,\frac{d-s}{d-c}c+\frac{s-c}{d-c}d\right) \right\vert ^{q} \\
&&+\frac{t-a}{b-a}\left\vert \frac{\partial ^{2}f}{\partial t\partial s}%
\left( b,\frac{d-s}{d-c}c+\frac{s-c}{d-c}d\right) \right\vert ^{q}
\end{eqnarray*}%
and%
\begin{eqnarray}
&&\left\vert \frac{\partial ^{2}f}{\partial t\partial s}\left( \frac{b-t}{b-a%
}a+\frac{t-a}{b-a}b,\frac{d-s}{d-c}c+\frac{s-c}{d-c}d\right) \right\vert ^{q}
\label{2.4} \\
&\leq &\left( \frac{b-t}{b-a}\right) \left( \frac{d-s}{d-c}\right)
\left\vert \frac{\partial ^{2}f}{\partial t\partial s}\left( a,c\right)
\right\vert ^{q}  \notag \\
&&+\left( \frac{b-t}{b-a}\right) \left( \frac{s-c}{d-c}\right) \left\vert 
\frac{\partial ^{2}f}{\partial t\partial s}\left( a,d\right) \right\vert ^{q}
\notag \\
&&+\left( \frac{t-a}{b-a}\right) \left( \frac{d-s}{d-c}\right) \left\vert 
\frac{\partial ^{2}f}{\partial t\partial s}\left( b,c\right) \right\vert ^{q}
\notag \\
&&+\left( \frac{t-a}{b-a}\right) \left( \frac{s-c}{d-c}\right) \left\vert 
\frac{\partial ^{2}f}{\partial t\partial s}\left( b,d\right) \right\vert
^{q}.  \notag
\end{eqnarray}%
Using inequality of (\ref{2.4}) in (\ref{2.3}), we get%
\begin{eqnarray*}
&&\left\vert f\left( \frac{a+b}{2},\frac{c+d}{2}\right) +\frac{1}{\left(
b-a\right) \left( d-c\right) }\int_{a}^{b}\int_{c}^{d}f\left( x,y\right)
dydx\right. \\
&&\left. -\frac{1}{\left( d-c\right) }\int_{c}^{d}f\left( \frac{a+b}{2}%
,y\right) dy-\frac{1}{\left( b-a\right) }\int_{a}^{b}f\left( x,\frac{c+d}{2}%
\right) dx\right\vert \\
&\leq &\frac{\left( b-a\right) \left( d-c\right) }{4\left( p+1\right) ^{%
\frac{2}{p}}}\left( \frac{\left\vert \frac{\partial ^{2}f}{\partial
t\partial s}\left( a,c\right) \right\vert ^{q}+\left\vert \frac{\partial
^{2}f}{\partial t\partial s}\left( b,c\right) \right\vert ^{q}+\left\vert 
\frac{\partial ^{2}f}{\partial t\partial s}\left( a,d\right) \right\vert
^{q}+\left\vert \frac{\partial ^{2}f}{\partial t\partial s}\left( b,d\right)
\right\vert ^{q}}{4}\right) ^{\frac{1}{q}}
\end{eqnarray*}%
where we have used the fact that%
\begin{equation*}
\left( \int_{a}^{b}\int_{c}^{d}\left\vert p\left( x,t\right) q\left(
y,s\right) \right\vert ^{p}dtds\right) ^{\frac{1}{p}}=\frac{\left[ \left(
b-a\right) \left( d-c\right) \right] ^{1+\frac{1}{p}}}{4\left( p+1\right) ^{%
\frac{2}{p}}}.
\end{equation*}%
This completes the proof.
\end{proof}

\begin{remark}
Suppose that all the assumptions of Theorem \ref{t.2.2} are satisfied. If we
choose $\frac{\partial ^{2}f}{\partial t\partial s}$ is bounded, i.e.,%
\begin{equation*}
\left\Vert \frac{\partial ^{2}f\left( t,s\right) }{\partial t\partial s}%
\right\Vert _{\infty }=\underset{\left( t,s\right) \in \left( a,b\right)
\times \left( c,d\right) }{\sup }\left\vert \frac{\partial ^{2}f\left(
t,s\right) }{\partial t\partial s}\right\vert <\infty ,
\end{equation*}%
we get%
\begin{eqnarray}
&&\left\vert f\left( \frac{a+b}{2},\frac{c+d}{2}\right) \right.
\label{r.2.2} \\
&&-\frac{1}{\left( d-c\right) }\int_{c}^{d}f\left( \frac{a+b}{2},y\right) dy-%
\frac{1}{\left( b-a\right) }\int_{a}^{b}f\left( x,\frac{c+d}{2}\right) dx 
\notag \\
&&\left. +\frac{1}{\left( b-a\right) \left( d-c\right) }\int_{a}^{b}%
\int_{c}^{d}f\left( x,y\right) dydx\right\vert  \notag \\
&\leq &\frac{\left( b-a\right) \left( d-c\right) }{4\left( p+1\right) ^{%
\frac{2}{p}}}\left\Vert \frac{\partial ^{2}f\left( t,s\right) }{\partial
t\partial s}\right\Vert _{\infty }.  \notag
\end{eqnarray}
\end{remark}

\begin{theorem}
\label{t.2.3} Let $f:\Delta =\left[ a,b\right] \times \left[ c,d\right]
\rightarrow 
\mathbb{R}
$ be a partial differentiable mapping on $\Delta =\left[ a,b\right] \times %
\left[ c,d\right] .$ If $\left\vert \frac{\partial ^{2}f}{\partial t\partial
s}\right\vert ^{q},$ $q\geq 1,$ is a convex function on the co-ordinates on $%
\Delta ,$ then the following inequality holds;%
\begin{eqnarray}
&&\left\vert f\left( \frac{a+b}{2},\frac{c+d}{2}\right) +\frac{1}{\left(
b-a\right) \left( d-c\right) }\int_{a}^{b}\int_{c}^{d}f\left( x,y\right)
dydx\right.  \label{2.5} \\
&&\left. -\frac{1}{\left( d-c\right) }\int_{c}^{d}f\left( \frac{a+b}{2}%
,y\right) dy-\frac{1}{\left( b-a\right) }\int_{a}^{b}f\left( x,\frac{c+d}{2}%
\right) dx\right\vert  \notag \\
&\leq &\frac{\left( b-a\right) \left( d-c\right) }{16}  \notag \\
&&\times \left( \frac{\left\vert \frac{\partial ^{2}f}{\partial t\partial s}%
\left( a,c\right) \right\vert ^{q}+\left\vert \frac{\partial ^{2}f}{\partial
t\partial s}\left( b,c\right) \right\vert ^{q}+\left\vert \frac{\partial
^{2}f}{\partial t\partial s}\left( a,d\right) \right\vert ^{q}+\left\vert 
\frac{\partial ^{2}f}{\partial t\partial s}\left( b,d\right) \right\vert ^{q}%
}{4}\right) ^{\frac{1}{q}}.  \notag
\end{eqnarray}
\end{theorem}

\begin{proof}
From Lemma 2, we have%
\begin{eqnarray*}
&&\left\vert f\left( \frac{a+b}{2},\frac{c+d}{2}\right) +\frac{1}{\left(
b-a\right) \left( d-c\right) }\int_{a}^{b}\int_{c}^{d}f\left( x,y\right)
dydx\right. \\
&&\left. -\frac{1}{\left( d-c\right) }\int_{c}^{d}f\left( \frac{a+b}{2}%
,y\right) dy-\frac{1}{\left( b-a\right) }\int_{a}^{b}f\left( x,\frac{c+d}{2}%
\right) dx\right\vert \\
&\leq &\frac{1}{\left( b-a\right) \left( d-c\right) } \\
&&\int_{a}^{b}\int_{c}^{d}\left\vert p\left( x,t\right) q\left( y,s\right)
\right\vert \left\vert \frac{\partial ^{2}f}{\partial t\partial s}\left( 
\frac{b-t}{b-a}a+\frac{t-a}{b-a}b,\frac{d-s}{d-c}c+\frac{s-c}{d-c}d\right)
\right\vert dsdt
\end{eqnarray*}%
By applying the well-known Power mean inequality for double integrals, then
one has%
\begin{eqnarray}
&&\left\vert f\left( \frac{a+b}{2},\frac{c+d}{2}\right) +\frac{1}{\left(
b-a\right) \left( d-c\right) }\int_{a}^{b}\int_{c}^{d}f\left( x,y\right)
dydx\right.  \notag \\
&&\left. -\frac{1}{\left( d-c\right) }\int_{c}^{d}f\left( \frac{a+b}{2}%
,y\right) dy-\frac{1}{\left( b-a\right) }\int_{a}^{b}f\left( x,\frac{c+d}{2}%
\right) dx\right\vert  \notag \\
&\leq &\frac{1}{\left( b-a\right) \left( d-c\right) }\times \left(
\int_{a}^{b}\int_{c}^{d}\left\vert p\left( x,t\right) q\left( y,s\right)
\right\vert dsdt\right) ^{1-\frac{1}{q}}  \label{z} \\
&&\left( \int\limits_{a}^{b}\int\limits_{c}^{d}\left\vert p\left( x,t\right)
q\left( y,s\right) \right\vert \left\vert \frac{\partial ^{2}f}{\partial
t\partial s}\left( \frac{b-t}{b-a}a+\frac{t-a}{b-a}b,\frac{d-s}{d-c}c+\frac{%
s-c}{d-c}d\right) \right\vert ^{q}dsdt\right) ^{\frac{1}{q}}.  \notag
\end{eqnarray}%
Since $\left\vert \frac{\partial ^{2}f}{\partial t\partial s}\right\vert
^{q} $ is a co-ordinated convex function on $\Delta ,$ we can write for all $%
\left( t,s\right) \in \left[ a,b\right] \times \left[ c,d\right] $%
\begin{eqnarray}
&&\left\vert \frac{\partial ^{2}f}{\partial t\partial s}\left( \frac{b-t}{b-a%
}a+\frac{t-a}{b-a}b,\frac{d-s}{d-c}c+\frac{s-c}{d-c}d\right) \right\vert ^{q}
\label{2.7} \\
&\leq &\left( \frac{b-t}{b-a}\right) \left( \frac{d-s}{d-c}\right)
\left\vert \frac{\partial ^{2}f}{\partial t\partial s}\left( a,c\right)
\right\vert ^{q}  \notag \\
&&+\left( \frac{b-t}{b-a}\right) \left( \frac{s-c}{d-c}\right) \left\vert 
\frac{\partial ^{2}f}{\partial t\partial s}\left( a,d\right) \right\vert ^{q}
\notag \\
&&+\left( \frac{t-a}{b-a}\right) \left( \frac{d-s}{d-c}\right) \left\vert 
\frac{\partial ^{2}f}{\partial t\partial s}\left( b,c\right) \right\vert ^{q}
\notag \\
&&+\left( \frac{t-a}{b-a}\right) \left( \frac{s-c}{d-c}\right) \left\vert 
\frac{\partial ^{2}f}{\partial t\partial s}\left( b,d\right) \right\vert
^{q}.  \notag
\end{eqnarray}

If we use ($\ref{2.7})$ in $\left( \ref{z}\right) $, we get%
\begin{eqnarray*}
&&\left\vert f\left( \frac{a+b}{2},\frac{c+d}{2}\right) +\frac{1}{\left(
b-a\right) \left( d-c\right) }\int_{a}^{b}\int_{c}^{d}f\left( x,y\right)
dydx\right. \\
&&\left. -\frac{1}{\left( d-c\right) }\int_{c}^{d}f\left( \frac{a+b}{2}%
,y\right) dy-\frac{1}{\left( b-a\right) }\int_{a}^{b}f\left( x,\frac{c+d}{2}%
\right) dx\right\vert \\
&\leq &\frac{1}{\left( b-a\right) \left( d-c\right) }\left\{ \left(
\int_{a}^{b}\int_{c}^{d}\left\vert p\left( x,t\right) q\left( y,s\right)
\right\vert dsdt\right) ^{1-\frac{1}{q}}\right. \\
&&\times \left( \int_{a}^{b}\int_{c}^{d}\left\vert p\left( x,t\right)
q\left( y,s\right) \right\vert \left\vert \frac{\partial ^{2}f}{\partial
t\partial s}\left( \frac{b-t}{b-a}a+\frac{t-a}{b-a}b,\frac{d-s}{d-c}c+\frac{%
s-c}{d-c}d\right) \right\vert ^{q}dsdt\right) ^{\frac{1}{q}}
\end{eqnarray*}%
hence it follows that%
\begin{eqnarray*}
&&\left\vert f\left( \frac{a+b}{2},\frac{c+d}{2}\right) +\frac{1}{\left(
b-a\right) \left( d-c\right) }\int_{a}^{b}\int_{c}^{d}f\left( x,y\right)
dydx\right. \\
&&\left. -\frac{1}{\left( d-c\right) }\int_{c}^{d}f\left( \frac{a+b}{2}%
,y\right) dy-\frac{1}{\left( b-a\right) }\int_{a}^{b}f\left( x,\frac{c+d}{2}%
\right) dx\right\vert \\
&\leq &\frac{1}{\left( b-a\right) \left( d-c\right) }\left\{ \left(
\int_{a}^{b}\int_{c}^{d}\left\vert p\left( x,t\right) q\left( y,s\right)
\right\vert dsdt\right) ^{1-\frac{1}{q}}\right. \\
&&\times \left( \int_{a}^{b}\int_{c}^{d}\left\vert p\left( x,t\right)
q\left( y,s\right) \right\vert \left[ \left( \frac{b-t}{b-a}\right) \left( 
\frac{d-s}{d-c}\right) \left\vert \frac{\partial ^{2}f}{\partial t\partial s}%
\left( a,c\right) \right\vert ^{q}+\left( \frac{b-t}{b-a}\right) \left( 
\frac{s-c}{d-c}\right) \left\vert \frac{\partial ^{2}f}{\partial t\partial s}%
\left( a,d\right) \right\vert ^{q}\right. \right. \\
&&\left. \left. +\left( \frac{t-a}{b-a}\right) \left( \frac{d-s}{d-c}\right)
\left\vert \frac{\partial ^{2}f}{\partial t\partial s}\left( b,c\right)
\right\vert ^{q}+\left( \frac{t-a}{b-a}\right) \left( \frac{s-c}{d-c}\right)
\left\vert \frac{\partial ^{2}f}{\partial t\partial s}\left( b,d\right)
\right\vert ^{q}\right] \right) ^{\frac{1}{q}}.
\end{eqnarray*}%
Computing the above integrals and using the fact that%
\begin{equation*}
\left( \int_{a}^{b}\int_{c}^{d}\left\vert p\left( x,t\right) q\left(
y,s\right) \right\vert dtds\right) ^{1-\frac{1}{q}}=\left( \frac{\left(
b-a\right) ^{2}\left( d-c\right) ^{2}}{16}\right) ^{1-\frac{1}{q}}.
\end{equation*}%
This completes the proof.
\end{proof}

\section{INEQUALITIES\ FOR CO-ORDINATED $s-$CONVEX FUNCTIONS}

To prove our main results we need the following lemma:

\begin{lemma}
Let $f:\Delta \subset 
\mathbb{R}
^{2}\rightarrow 
\mathbb{R}
$ be an absolutely continuous function on $\Delta $ where $a<b,$ $c<d$ and $%
t,\lambda \in \left[ 0,1\right] $, if $\frac{\partial ^{2}f}{\partial
t\partial \lambda }\in L\left( \Delta \right) $, then the following equality
holds:%
\begin{eqnarray*}
&&\frac{f\left( a,c\right) +r_{2}f\left( a,d\right) +r_{1}f\left( b,c\right)
+r_{1}r_{2}f\left( b,d\right) }{\left( r_{1}+1\right) \left( r_{2}+1\right) }
\\
&&+\frac{1}{(b-a)(d-c)}\int_{a}^{b}\int_{c}^{d}f(x,y)dxdy \\
&&-\left( \frac{r_{2}}{r_{2}+1}\right) \frac{1}{d-c}\int_{c}^{d}f(b,y)dy-%
\left( \frac{1}{r_{1}+1}\right) \frac{1}{d-c}\int_{c}^{d}f(a,y)dy \\
&&-\left( \frac{r_{2}}{r_{2}+1}\right) \frac{1}{b-a}\int_{a}^{b}f(x,d)dx-%
\left( \frac{1}{r_{2}+1}\right) \frac{1}{b-a}\int_{a}^{b}f(x,c)dx \\
&=&\frac{(b-a)(d-c)}{\left( r_{1}+1\right) \left( r_{2}+1\right) } \\
&&\times \int_{0}^{1}\int_{0}^{1}\left( \left( r_{1}+1\right) t-1\right)
\left( \left( r_{2}+1\right) \lambda -1\right) \frac{\partial ^{2}f}{%
\partial t\partial \lambda }\left( tb+\left( 1-t\right) a,\lambda d+\left(
1-\lambda \right) c\right) dtd\lambda
\end{eqnarray*}%
for some fixed $r_{1},r_{2}\in \left[ 0,1\right] .$
\end{lemma}

\begin{proof}
Integration by parts, we get%
\begin{eqnarray*}
&&\int_{0}^{1}\int_{0}^{1}\left( \left( r_{1}+1\right) t-1\right) \left(
\left( r_{2}+1\right) \lambda -1\right) \frac{\partial ^{2}f}{\partial
t\partial \lambda }\left( tb+\left( 1-t\right) a,\lambda d+\left( 1-\lambda
\right) c\right) dtd\lambda \\
&=&\int_{0}^{1}\left( \left( r_{2}+1\right) \lambda -1\right) \left[
\int_{0}^{1}\left( \left( r_{1}+1\right) t-1\right) \frac{\partial ^{2}f}{%
\partial t\partial \lambda }\left( tb+\left( 1-t\right) a,\lambda d+\left(
1-\lambda \right) c\right) dtd\lambda \right] \\
&=&\int_{0}^{1}\left( \left( r_{2}+1\right) \lambda -1\right) \left[ \frac{%
\left( \left( r_{1}+1\right) t-1\right) }{\left( b-a\right) }\frac{\partial f%
}{\partial \lambda }\left( tb+\left( 1-t\right) a,\lambda d+\left( 1-\lambda
\right) c\right) \right\vert _{0}^{1} \\
&&\left. -\frac{r_{1}+1}{b-a}\int_{0}^{1}\frac{\partial f}{\partial \lambda }%
\left( tb+\left( 1-t\right) a,\lambda d+\left( 1-\lambda \right) c\right) dt%
\right] d\lambda
\end{eqnarray*}%
\begin{eqnarray*}
&&\int_{0}^{1}\int_{0}^{1}\left( \left( r_{1}+1\right) t-1\right) \left(
\left( r_{2}+1\right) \lambda -1\right) \frac{\partial ^{2}f}{\partial
t\partial \lambda }\left( tb+\left( 1-t\right) a,\lambda d+\left( 1-\lambda
\right) c\right) dtd\lambda \\
&=&\int_{0}^{1}\left( \left( r_{2}+1\right) \lambda -1\right) \left[ \frac{%
r_{1}}{b-a}\frac{\partial f}{\partial \lambda }\left( b,\lambda d+\left(
1-\lambda \right) c\right) +\frac{1}{b-a}\frac{\partial f}{\partial \lambda }%
\left( a,\lambda d+\left( 1-\lambda \right) c\right) \right. \\
&&\left. -\frac{r_{1}+1}{b-a}\int_{0}^{1}\frac{\partial f}{\partial \lambda }%
\left( tb+\left( 1-t\right) a,\lambda d+\left( 1-\lambda \right) c\right) dt%
\right] d\lambda
\end{eqnarray*}%
Again by integration by parts, we have%
\begin{eqnarray*}
&&\int_{0}^{1}\left( \left( r_{2}+1\right) \lambda -1\right) \left[ \frac{%
r_{1}}{b-a}\frac{\partial f}{\partial \lambda }\left( b,\lambda d+\left(
1-\lambda \right) c\right) +\frac{1}{b-a}\frac{\partial f}{\partial \lambda }%
\left( a,\lambda d+\left( 1-\lambda \right) c\right) \right. \\
&&\left. -\frac{r_{1}+1}{b-a}\int_{0}^{1}\frac{\partial f}{\partial \lambda }%
\left( tb+\left( 1-t\right) a,\lambda d+\left( 1-\lambda \right) c\right) dt%
\right] d\lambda \\
&=&\left. \frac{r_{1}}{b-a}\frac{\left( \left( r_{2}+1\right) \lambda
-1\right) }{d-c}f\left( b,\lambda d+\left( 1-\lambda \right) c\right)
\right\vert _{0}^{1}-\frac{r_{1}\left( r_{2}+1\right) }{(b-a)(d-c)}%
\int_{0}^{1}f\left( b,\lambda d+\left( 1-\lambda \right) c\right) d\lambda \\
&&\left. +\frac{1}{b-a}\frac{\left( \left( r_{2}+1\right) \lambda -1\right) 
}{d-c}f\left( a,\lambda d+\left( 1-\lambda \right) c\right) \right\vert
_{0}^{1}-\frac{\left( r_{2}+1\right) }{(b-a)(d-c)}\int_{0}^{1}f\left(
a,\lambda d+\left( 1-\lambda \right) c\right) d\lambda \\
&&-\frac{r_{1}+1}{b-a}\int_{0}^{1}\left[ \int_{0}^{1}\left( \left(
r_{2}+1\right) \lambda -1\right) \frac{\partial f}{\partial \lambda }\left(
tb+\left( 1-t\right) a,\lambda d+\left( 1-\lambda \right) c\right) d\lambda %
\right] dt.
\end{eqnarray*}%
Computing these integrals and by using the results, we obtain%
\begin{eqnarray*}
&&\int_{0}^{1}\int_{0}^{1}\left( \left( r_{1}+1\right) t-1\right) \left(
\left( r_{2}+1\right) \lambda -1\right) \frac{\partial ^{2}f}{\partial
t\partial \lambda }\left( tb+\left( 1-t\right) a,\lambda d+\left( 1-\lambda
\right) c\right) dtd\lambda \\
&=&\frac{1}{(b-a)(d-c)}\left[ f\left( a,c\right) +r_{2}f\left( a,d\right)
+r_{1}f\left( b,c\right) +r_{1}r_{2}f\left( b,d\right) \right. \\
&&-r_{1}\left( r_{2}+1\right) \int_{0}^{1}f(b,\lambda d+\left( 1-\lambda
\right) c)d\lambda -\left( r_{2}+1\right) \int_{0}^{1}f(a,\lambda d+\left(
1-\lambda \right) c)d\lambda \\
&&-r_{2}\left( r_{1}+1\right) \int_{0}^{1}f(tb+\left( 1-t\right)
a,d)dt-\left( r_{2}+1\right) \int_{0}^{1}f(tb+\left( 1-t\right) a,c)dt \\
&&\left. \left( r_{1}+1\right) \left( r_{2}+1\right)
\int_{0}^{1}\int_{0}^{1}f\left( tb+\left( 1-t\right) a,\lambda d+\left(
1-\lambda \right) c\right) dtd\lambda \right] .
\end{eqnarray*}%
Using the change of the variable $x=tb+\left( 1-t\right) a$ and $y=\lambda
d+\left( 1-\lambda \right) c$ for $t,\lambda \in \left[ 0,1\right] $ and
multiplying the both sides by $\frac{(b-a)(d-c)}{\left( r_{1}+1\right)
\left( r_{2}+1\right) },$ we get the required result.
\end{proof}

\begin{theorem}
Let $f:\Delta =[a,b]\times \lbrack c,d]\subset \lbrack 0,\infty
)^{2}\rightarrow \lbrack 0,\infty )$ be an absolutely continuous function on 
$\Delta $. If $\left\vert \frac{\partial ^{2}f}{\partial t\partial \lambda }%
\right\vert $ is $s-$convex function on the co-ordinates on $\Delta ,$ then
one has the inequality:%
\begin{eqnarray}
&&\left\vert \frac{f\left( a,c\right) +r_{2}f\left( a,d\right) +r_{1}f\left(
b,c\right) +r_{1}r_{2}f\left( b,d\right) }{\left( r_{1}+1\right) \left(
r_{2}+1\right) }\right.  \notag \\
&&-\frac{1}{d-c}\left[ \left( \frac{r_{2}}{r_{2}+1}\right)
\int_{c}^{d}f(b,y)dy+\left( \frac{1}{r_{1}+1}\right) \int_{c}^{d}f(a,y)dy%
\right]  \notag \\
&&-\frac{1}{b-a}\left[ \left( \frac{r_{2}}{r_{2}+1}\right)
\int_{a}^{b}f(x,d)dx+\left( \frac{1}{r_{2}+1}\right) \int_{a}^{b}f(x,c)dx%
\right]  \label{9} \\
&&\left. +\frac{1}{(b-a)(d-c)}\int_{a}^{b}\int_{c}^{d}f(x,y)dxdy\right\vert 
\notag \\
&\leq &\frac{(b-a)(d-c)}{\left( r_{1}+1\right) \left( r_{2}+1\right) \left(
s+1\right) ^{2}\left( s+2\right) ^{2}}  \notag \\
&&\times \left( MN\left\vert \frac{\partial ^{2}f}{\partial t\partial
\lambda }\right\vert (a,c)+LN\left\vert \frac{\partial ^{2}f}{\partial
t\partial \lambda }\right\vert (a,d)\right.  \notag \\
&&\left. +KM\left\vert \frac{\partial ^{2}f}{\partial t\partial \lambda }%
\right\vert (b,c)+KL\left\vert \frac{\partial ^{2}f}{\partial t\partial
\lambda }\right\vert (b,d)\right)  \notag
\end{eqnarray}%
where%
\begin{eqnarray*}
M &=&\left( s+1+2\left( r_{1}+1\right) \left( \frac{r_{1}}{r_{1}+1}\right)
^{s+2}-r_{1}\right) \\
N &=&\left( s+1+2\left( r_{2}+1\right) \left( \frac{r_{2}}{r_{2}+1}\right)
^{s+2}-r_{2}\right) \\
L &=&\left( r_{2}\left( s+1\right) +2\left( \frac{1}{r_{2}+1}\right)
^{s+1}-1\right) \\
K &=&\left( r_{1}\left( s+1\right) +2\left( \frac{1}{r_{1}+1}\right)
^{s+1}-1\right) .
\end{eqnarray*}
\end{theorem}

\begin{proof}
From Lemma 3, we can write;%
\begin{eqnarray*}
&&\left\vert \frac{f\left( a,c\right) +r_{2}f\left( a,d\right) +r_{1}f\left(
b,c\right) +r_{1}r_{2}f\left( b,d\right) }{\left( r_{1}+1\right) \left(
r_{2}+1\right) }+\frac{1}{(b-a)(d-c)}\int_{a}^{b}\int_{c}^{d}f(x,y)dxdy%
\right. \\
&&-\left( \frac{r_{2}}{r_{2}+1}\right) \frac{1}{d-c}\int_{c}^{d}f(b,y)dy-%
\left( \frac{1}{r_{1}+1}\right) \frac{1}{d-c}\int_{c}^{d}f(a,y)dy \\
&&\left. -\left( \frac{r_{2}}{r_{2}+1}\right) \frac{1}{b-a}%
\int_{a}^{b}f(x,d)dx-\left( \frac{1}{r_{2}+1}\right) \frac{1}{b-a}%
\int_{a}^{b}f(x,c)dx\right\vert \\
&=&\frac{(b-a)(d-c)}{\left( r_{1}+1\right) \left( r_{2}+1\right) } \\
&&\times \int_{0}^{1}\int_{0}^{1}\left\vert \left( \left( r_{1}+1\right)
t-1\right) \left( \left( r_{2}+1\right) \lambda -1\right) \right\vert
\left\vert \frac{\partial ^{2}f}{\partial t\partial \lambda }\left(
tb+\left( 1-t\right) a,\lambda d+\left( 1-\lambda \right) c\right)
\right\vert dtd\lambda
\end{eqnarray*}%
By using co-ordinated $s-$convexity of $f,$ we have%
\begin{eqnarray*}
&&\left\vert \frac{f\left( a,c\right) +r_{2}f\left( a,d\right) +r_{1}f\left(
b,c\right) +r_{1}r_{2}f\left( b,d\right) }{\left( r_{1}+1\right) \left(
r_{2}+1\right) }\right. \\
&&-\left( \frac{r_{2}}{r_{2}+1}\right) \frac{1}{d-c}\int_{c}^{d}f(b,y)dy-%
\left( \frac{1}{r_{1}+1}\right) \frac{1}{d-c}\int_{c}^{d}f(a,y)dy \\
&&-\left( \frac{r_{2}}{r_{2}+1}\right) \frac{1}{b-a}\int_{a}^{b}f(x,d)dx-%
\left( \frac{1}{r_{2}+1}\right) \frac{1}{b-a}\int_{a}^{b}f(x,c)dx \\
&&\left. +\frac{1}{(b-a)(d-c)}\int_{a}^{b}\int_{c}^{d}f(x,y)dxdy\right\vert
\\
&=&\frac{(b-a)(d-c)}{\left( r_{1}+1\right) \left( r_{2}+1\right) } \\
&&\times \int_{0}^{1}\left[ \int_{0}^{1}\left\vert \left( \left(
r_{1}+1\right) t-1\right) \left( \left( r_{2}+1\right) \lambda -1\right)
\right\vert \right. \\
&&\left. \left\{ t^{s}\left\vert \frac{\partial ^{2}f}{\partial t\partial
\lambda }\left( b,\lambda d+\left( 1-\lambda \right) c\right) \right\vert
+\left( 1-t\right) ^{s}\left\vert \frac{\partial ^{2}f}{\partial t\partial
\lambda }\left( a,\lambda d+\left( 1-\lambda \right) c\right) \right\vert
\right\} dt\right] d\lambda .
\end{eqnarray*}%
By calculating the above integrals, we get%
\begin{eqnarray}
&&\int_{0}^{1}\left\vert \left( \left( r_{1}+1\right) t-1\right) \right\vert
\left\{ t^{s}\left\vert \frac{\partial ^{2}f}{\partial t\partial \lambda }%
\left( b,\lambda d+\left( 1-\lambda \right) c\right) \right\vert \right.
\label{3} \\
&&\left. +\left( 1-t\right) ^{s}\left\vert \frac{\partial ^{2}f}{\partial
t\partial \lambda }\left( a,\lambda d+\left( 1-\lambda \right) c\right)
\right\vert \right\} dt  \notag \\
&=&\int_{0}^{\frac{1}{r_{1}+1}}\left( 1-\left( r_{1}+1\right) t\right)
\left\{ t^{s}\left\vert \frac{\partial ^{2}f}{\partial t\partial \lambda }%
\left( b,\lambda d+\left( 1-\lambda \right) c\right) \right\vert \right. 
\notag \\
&&\left. +\left( 1-t\right) ^{s}\left\vert \frac{\partial ^{2}f}{\partial
t\partial \lambda }\left( a,\lambda d+\left( 1-\lambda \right) c\right)
\right\vert \right\} dt  \notag \\
&&+\int_{\frac{1}{r_{1}+1}}^{1}\left( \left( r_{1}+1\right) t-1\right)
\left\{ t^{s}\left\vert \frac{\partial ^{2}f}{\partial t\partial \lambda }%
\left( b,\lambda d+\left( 1-\lambda \right) c\right) \right\vert \right. 
\notag \\
&&\left. +\left( 1-t\right) ^{s}\left\vert \frac{\partial ^{2}f}{\partial
t\partial \lambda }\left( a,\lambda d+\left( 1-\lambda \right) c\right)
\right\vert \right\} dt  \notag \\
&=&\frac{1}{\left( s+1\right) \left( s+2\right) }\left[ \left( r_{1}\left(
s+1\right) +2\left( \frac{1}{r_{1}+1}\right) ^{s+1}-1\right) \left\vert 
\frac{\partial ^{2}f}{\partial t\partial \lambda }\left( b,\lambda d+\left(
1-\lambda \right) c\right) \right\vert \right.  \notag \\
&&\left. +\left( s+1+2\left( r_{1}+1\right) \left( \frac{r_{1}}{r_{1}+1}%
\right) ^{s+2}-r_{1}\right) \left\vert \frac{\partial ^{2}f}{\partial
t\partial \lambda }\left( a,\lambda d+\left( 1-\lambda \right) c\right)
\right\vert \right] .  \notag
\end{eqnarray}%
By a similar argument for other integrals, by using co-ordinated $s-$%
convexity of $f,$ we get%
\begin{eqnarray*}
&&\int_{0}^{1}\left\vert \left( \left( r_{2}+1\right) t-1\right) \right\vert
\left\{ \left\vert \frac{\partial ^{2}f}{\partial t\partial \lambda }\left(
b,\lambda d+\left( 1-\lambda \right) c\right) \right\vert +\left\vert \frac{%
\partial ^{2}f}{\partial t\partial \lambda }\left( a,\lambda d+\left(
1-\lambda \right) c\right) \right\vert \right\} d\lambda \\
&=&\int_{0}^{\frac{1}{r_{2}+1}}\left( 1-\left( r_{2}+1\right) t\right)
\left\{ \lambda ^{s}\left\vert \frac{\partial ^{2}f}{\partial t\partial
\lambda }\left( b,d\right) \right\vert +\left( 1-\lambda \right)
^{s}\left\vert \frac{\partial ^{2}f}{\partial t\partial \lambda }\left(
b,c\right) \right\vert \right\} d\lambda \\
&&+\int_{\frac{1}{r_{2}+1}}^{1}\left( \left( r_{2}+1\right) t-1\right)
\left\{ \lambda ^{s}\left\vert \frac{\partial ^{2}f}{\partial t\partial
\lambda }\left( a,d\right) \right\vert +\left( 1-\lambda \right)
^{s}\left\vert \frac{\partial ^{2}f}{\partial t\partial \lambda }\left(
a,c\right) \right\vert \right\} d\lambda \\
&=&\frac{1}{\left( s+1\right) \left( s+2\right) }\left[ r_{2}\left(
s+1\right) +2\left( \frac{1}{r_{2}+1}\right) ^{s+1}-1\right] \left\vert 
\frac{\partial ^{2}f}{\partial t\partial \lambda }\left( b,d\right)
\right\vert \\
&&+\frac{1}{\left( s+1\right) \left( s+2\right) }\left[ r_{2}\left(
s+1\right) +2\left( \frac{1}{r_{2}+1}\right) ^{s+1}-1\right] \left\vert 
\frac{\partial ^{2}f}{\partial t\partial \lambda }\left( a,d\right)
\right\vert \\
&&+\frac{1}{\left( s+1\right) \left( s+2\right) }\left[ s+1+2\left(
r_{2}+1\right) \left( \frac{r_{2}}{r_{2}+1}\right) ^{s+2}-r_{1}\right]
\left\vert \frac{\partial ^{2}f}{\partial t\partial \lambda }\left(
b,c\right) \right\vert \\
&&+\frac{1}{\left( s+1\right) \left( s+2\right) }\left[ s+1+2\left(
r_{1}+1\right) \left( \frac{r_{2}}{r_{2}+1}\right) ^{s+2}-r_{1}\right]
\left\vert \frac{\partial ^{2}f}{\partial t\partial \lambda }\left(
a,c\right) \right\vert .
\end{eqnarray*}%
By using these in (\ref{3}), we obtain the inequality (\ref{9}).
\end{proof}

\begin{corollary}
(1) If we choose $r_{1}=r_{2}=1$ in (\ref{9}), we have%
\begin{eqnarray}
&&\left\vert \frac{f\left( a,c\right) +f\left( a,d\right) +f\left(
b,c\right) +f\left( b,d\right) }{4}\right.  \label{4} \\
&&-\frac{1}{2}\left[ \frac{1}{d-c}\int_{c}^{d}\left[ f(b,y)+f(a,y)\right] dy%
\right] -\frac{1}{2}\left[ \frac{1}{b-a}\int_{a}^{b}\left[ f(x,d)+f(x,c)%
\right] dx\right]  \notag \\
&&\left. +\frac{1}{(b-a)(d-c)}\int_{a}^{b}\int_{c}^{d}f(x,y)dxdy\right\vert 
\notag \\
&\leq &\frac{(b-a)(d-c)}{\left( s+1\right) ^{2}\left( s+2\right) ^{2}}\left(
s+\frac{1}{2^{s}}\right) ^{2}  \notag \\
&&\left( \left\vert \frac{\partial ^{2}f}{\partial t\partial \lambda }%
\right\vert (a,c)+\left\vert \frac{\partial ^{2}f}{\partial t\partial
\lambda }\right\vert (a,d)+\left\vert \frac{\partial ^{2}f}{\partial
t\partial \lambda }\right\vert (b,c)+\left\vert \frac{\partial ^{2}f}{%
\partial t\partial \lambda }\right\vert (b,d)\right)  \notag
\end{eqnarray}%
(2) If we choose $r_{1}=r_{2}=0$ in (\ref{9}), we have%
\begin{eqnarray*}
&&\left\vert f\left( a,c\right) -\frac{1}{d-c}\int_{c}^{d}f(a,y)dy-\frac{1}{%
b-a}\int_{a}^{b}f(x,c)dx\right. \\
&&\left. +\frac{1}{(b-a)(d-c)}\int_{a}^{b}\int_{c}^{d}f(x,y)dxdy\right\vert
\\
&\leq &\frac{(b-a)(d-c)}{\left( s+1\right) ^{2}\left( s+2\right) ^{2}} \\
&&\left( \left( s+1\right) ^{2}\left\vert \frac{\partial ^{2}f}{\partial
t\partial \lambda }\right\vert (a,c)+\left( s+1\right) \left\vert \frac{%
\partial ^{2}f}{\partial t\partial \lambda }\right\vert (a,d)+\left(
s+1\right) \left\vert \frac{\partial ^{2}f}{\partial t\partial \lambda }%
\right\vert (b,c)+\left\vert \frac{\partial ^{2}f}{\partial t\partial
\lambda }\right\vert (b,d)\right)
\end{eqnarray*}
\end{corollary}

\begin{remark}
If we choose $s=1$ in (\ref{4}), we get an improvement for the inequality (%
\ref{1.3}).
\end{remark}

\begin{theorem}
Let $f:\Delta =[a,b]\times \lbrack c,d]\subset \lbrack 0,\infty
)^{2}\rightarrow \lbrack 0,\infty )$ be an absolutely continuous function on 
$\Delta $. If $\left\vert \frac{\partial ^{2}f}{\partial t\partial \lambda }%
\right\vert ^{\frac{p}{p-1}}$ is $s-$convex function on the co-ordinates on $%
\Delta ,$ for some fixed $s\in (0,1]$ and $p>1,$ then one has the inequality:%
\begin{eqnarray}
&&\left\vert \frac{f\left( a,c\right) +r_{2}f\left( a,d\right) +r_{1}f\left(
b,c\right) +r_{1}r_{2}f\left( b,d\right) }{\left( r_{1}+1\right) \left(
r_{2}+1\right) }\right.  \label{a} \\
&&-\left( \frac{r_{2}}{r_{2}+1}\right) \frac{1}{d-c}\int_{c}^{d}f(b,y)dy-%
\left( \frac{1}{r_{1}+1}\right) \frac{1}{d-c}\int_{c}^{d}f(a,y)dy  \notag \\
&&-\left( \frac{r_{2}}{r_{2}+1}\right) \frac{1}{b-a}\int_{a}^{b}f(x,d)dx-%
\left( \frac{1}{r_{2}+1}\right) \frac{1}{b-a}\int_{a}^{b}f(x,c)dx  \notag \\
&&\left. +\frac{1}{(b-a)(d-c)}\int_{a}^{b}\int_{c}^{d}f(x,y)dxdy\right\vert 
\notag \\
&=&\frac{(b-a)(d-c)}{\left( r_{1}+1\right) \left( r_{2}+1\right) }\left( 
\frac{\left( 1+r_{1}^{\frac{p+1}{p}}\right) \left( 1+r_{2}^{\frac{p+1}{p}%
}\right) }{\left( r_{1}+1\right) ^{\frac{1}{p}}\left( r_{2}+1\right) ^{\frac{%
1}{p}}\left( p+1\right) ^{\frac{2}{p}}}\right)  \notag \\
&&\times \left( \frac{\left\vert \frac{\partial ^{2}f}{\partial t\partial
\lambda }\right\vert ^{q}(a,c)+\left\vert \frac{\partial ^{2}f}{\partial
t\partial \lambda }\right\vert ^{q}(a,d)+\left\vert \frac{\partial ^{2}f}{%
\partial t\partial \lambda }\right\vert ^{q}(b,c)+\left\vert \frac{\partial
^{2}f}{\partial t\partial \lambda }\right\vert ^{q}(b,d)}{\left( s+1\right)
^{2}}\right) ^{\frac{1}{q}}  \notag
\end{eqnarray}%
for some fixed $r_{1},r_{2}\in \left[ 0,1\right] ,$ where $q=\frac{p}{p-1}.$
\end{theorem}

\begin{proof}
Let $p>1.$ From Lemma 3 and using the H\"{o}lder inequality for double
integrals, we can write%
\begin{eqnarray*}
&&\left\vert \frac{f\left( a,c\right) +r_{2}f\left( a,d\right) +r_{1}f\left(
b,c\right) +r_{1}r_{2}f\left( b,d\right) }{\left( r_{1}+1\right) \left(
r_{2}+1\right) }\right. \\
&&-\left( \frac{r_{2}}{r_{2}+1}\right) \frac{1}{d-c}\int_{c}^{d}f(b,y)dy-%
\left( \frac{1}{r_{1}+1}\right) \frac{1}{d-c}\int_{c}^{d}f(a,y)dy \\
&&-\left( \frac{r_{2}}{r_{2}+1}\right) \frac{1}{b-a}\int_{a}^{b}f(x,d)dx-%
\left( \frac{1}{r_{2}+1}\right) \frac{1}{b-a}\int_{a}^{b}f(x,c)dx \\
&&\left. +\frac{1}{(b-a)(d-c)}\int_{a}^{b}\int_{c}^{d}f(x,y)dxdy\right\vert
\\
&=&\frac{(b-a)(d-c)}{\left( r_{1}+1\right) \left( r_{2}+1\right) }\left(
\int_{0}^{1}\int_{0}^{1}\left\vert \left( \left( r_{1}+1\right) t-1\right)
\left( \left( r_{2}+1\right) \lambda -1\right) \right\vert ^{p}dtd\lambda
\right) ^{\frac{1}{p}} \\
&&\times \left( \int_{0}^{1}\int_{0}^{1}\left\vert \frac{\partial ^{2}f}{%
\partial t\partial \lambda }\left( tb+\left( 1-t\right) a,\lambda d+\left(
1-\lambda \right) c\right) \right\vert ^{q}dtd\lambda \right) ^{\frac{1}{q}}.
\end{eqnarray*}%
Since $\left\vert \frac{\partial ^{2}f}{\partial t\partial \lambda }%
\right\vert ^{q}$ is $s-$convex function on the co-ordinates on $\Delta ,$
we can write for $t,\lambda \in \left[ 0,1\right] $%
\begin{eqnarray*}
&&\left\vert \frac{\partial ^{2}f}{\partial t\partial \lambda }\left(
tb+\left( 1-t\right) a,\lambda d+\left( 1-\lambda \right) c\right)
\right\vert ^{q} \\
&\leq &t^{s}\left\vert \frac{\partial ^{2}f}{\partial t\partial \lambda }%
\left( b,\lambda d+\left( 1-\lambda \right) c\right) \right\vert ^{q}+\left(
1-t\right) ^{s}\left\vert \frac{\partial ^{2}f}{\partial t\partial \lambda }%
\left( a,\lambda d+\left( 1-\lambda \right) c\right) \right\vert ^{q}
\end{eqnarray*}%
and%
\begin{eqnarray*}
&&\left\vert \frac{\partial ^{2}f}{\partial t\partial \lambda }\left(
tb+\left( 1-t\right) a,\lambda d+\left( 1-\lambda \right) c\right)
\right\vert ^{q} \\
&\leq &t^{s}\lambda ^{s}\left\vert \frac{\partial ^{2}f}{\partial t\partial
\lambda }\right\vert ^{q}(b,d)+t^{s}\left( 1-\lambda \right) ^{s}\left\vert 
\frac{\partial ^{2}f}{\partial t\partial \lambda }\right\vert ^{q}(b,c) \\
&&+\lambda ^{s}\left( 1-t\right) ^{s}\left\vert \frac{\partial ^{2}f}{%
\partial t\partial \lambda }\right\vert ^{q}(a,d)+\left( 1-\lambda \right)
^{s}\left( 1-t\right) ^{s}\left\vert \frac{\partial ^{2}f}{\partial
t\partial \lambda }\right\vert ^{q}(a,c)
\end{eqnarray*}%
thus, we obtain%
\begin{eqnarray*}
&&\left\vert \frac{f\left( a,c\right) +r_{2}f\left( a,d\right) +r_{1}f\left(
b,c\right) +r_{1}r_{2}f\left( b,d\right) }{\left( r_{1}+1\right) \left(
r_{2}+1\right) }+\frac{1}{(b-a)(d-c)}\int_{a}^{b}\int_{c}^{d}f(x,y)dxdy%
\right. \\
&&-\left( \frac{r_{2}}{r_{2}+1}\right) \frac{1}{d-c}\int_{c}^{d}f(b,y)dy-%
\left( \frac{1}{r_{1}+1}\right) \frac{1}{d-c}\int_{c}^{d}f(a,y)dy \\
&&\left. -\left( \frac{r_{2}}{r_{2}+1}\right) \frac{1}{b-a}%
\int_{a}^{b}f(x,d)dx-\left( \frac{1}{r_{2}+1}\right) \frac{1}{b-a}%
\int_{a}^{b}f(x,c)dx\right\vert \\
&=&\frac{(b-a)(d-c)}{\left( r_{1}+1\right) \left( r_{2}+1\right) }\left( 
\frac{\left( 1+r_{1}^{\frac{p+1}{p}}\right) \left( 1+r_{2}^{\frac{p+1}{p}%
}\right) }{\left( r_{1}+1\right) ^{\frac{1}{p}}\left( r_{2}+1\right) ^{\frac{%
1}{p}}\left( p+1\right) ^{\frac{2}{p}}}\right) \\
&&\times \left( \frac{\left\vert \frac{\partial ^{2}f}{\partial t\partial
\lambda }\right\vert ^{q}(a,c)+\left\vert \frac{\partial ^{2}f}{\partial
t\partial \lambda }\right\vert ^{q}(a,d)+\left\vert \frac{\partial ^{2}f}{%
\partial t\partial \lambda }\right\vert ^{q}(b,c)+\left\vert \frac{\partial
^{2}f}{\partial t\partial \lambda }\right\vert ^{q}(b,d)}{\left( s+1\right)
^{2}}\right) ^{\frac{1}{q}}.
\end{eqnarray*}%
Which completes the proof of the inequality (\ref{a}).
\end{proof}

\begin{corollary}
(1) Under the assumptions of Theorem 12, if we choose $r_{1}=r_{2}=1$ in (%
\ref{a})$,$ we have%
\begin{eqnarray}
&&\left\vert \frac{f\left( a,c\right) +f\left( a,d\right) +f\left(
b,c\right) +f\left( b,d\right) }{4}\right.  \label{5} \\
&&-\frac{1}{2}\left[ \frac{1}{d-c}\int_{c}^{d}\left[ f(b,y)+f(a,y)\right] dy+%
\frac{1}{b-a}\int_{a}^{b}\left[ f(x,d)+f(x,c)\right] dx\right]  \notag \\
&&\left. +\frac{1}{(b-a)(d-c)}\int_{a}^{b}\int_{c}^{d}f(x,y)dxdy\right\vert 
\notag \\
&=&\frac{(b-a)(d-c)}{4\left( p+1\right) ^{\frac{2}{p}}}  \notag \\
&&\times \left( \frac{\left\vert \frac{\partial ^{2}f}{\partial t\partial
\lambda }\right\vert ^{q}(a,c)+\left\vert \frac{\partial ^{2}f}{\partial
t\partial \lambda }\right\vert ^{q}(a,d)+\left\vert \frac{\partial ^{2}f}{%
\partial t\partial \lambda }\right\vert ^{q}(b,c)+\left\vert \frac{\partial
^{2}f}{\partial t\partial \lambda }\right\vert ^{q}(b,d)}{\left( s+1\right)
^{2}}\right) ^{\frac{1}{q}}.  \notag
\end{eqnarray}%
(2) Under the assumptions of Theorem 12, if we choose $r_{1}=r_{2}=0$ in (%
\ref{a})$,$ we have%
\begin{eqnarray*}
&&\left\vert f\left( a,c\right) -\frac{1}{d-c}\int_{c}^{d}f(a,y)dy-\frac{1}{%
b-a}\int_{a}^{b}f(x,c)dx\right. \\
&&\left. +\frac{1}{(b-a)(d-c)}\int_{a}^{b}\int_{c}^{d}f(x,y)dxdy\right\vert
\\
&=&\frac{(b-a)(d-c)}{\left( p+1\right) ^{\frac{2}{p}}} \\
&&\times \left( \frac{\left\vert \frac{\partial ^{2}f}{\partial t\partial
\lambda }\right\vert ^{q}(a,c)+\left\vert \frac{\partial ^{2}f}{\partial
t\partial \lambda }\right\vert ^{q}(a,d)+\left\vert \frac{\partial ^{2}f}{%
\partial t\partial \lambda }\right\vert ^{q}(b,c)+\left\vert \frac{\partial
^{2}f}{\partial t\partial \lambda }\right\vert ^{q}(b,d)}{\left( s+1\right)
^{2}}\right) ^{\frac{1}{q}}.
\end{eqnarray*}
\end{corollary}

\begin{remark}
If we choose $s=1$ in (\ref{5})$,$ we obtain an improvement for the
inequality (\ref{1.4}).
\end{remark}

\begin{theorem}
Let $f:\Delta =[a,b]\times \lbrack c,d]\subset \lbrack 0,\infty
)^{2}\rightarrow \lbrack 0,\infty )$ be an absolutely continuous function on 
$\Delta $. If $\left\vert \frac{\partial ^{2}f}{\partial t\partial \lambda }%
\right\vert ^{q}$ is $s-$convex function on the co-ordinates on $\Delta ,$
for some fixed $s\in (0,1]$ and $q\geq 1,$ then one has the inequality:%
\begin{eqnarray*}
&&\left\vert \frac{f\left( a,c\right) +r_{2}f\left( a,d\right) +r_{1}f\left(
b,c\right) +r_{1}r_{2}f\left( b,d\right) }{\left( r_{1}+1\right) \left(
r_{2}+1\right) }+\frac{1}{(b-a)(d-c)}\int_{a}^{b}\int_{c}^{d}f(x,y)dxdy%
\right. \\
&&-\left( \frac{r_{2}}{r_{2}+1}\right) \frac{1}{d-c}\int_{c}^{d}f(b,y)dy-%
\left( \frac{1}{r_{1}+1}\right) \frac{1}{d-c}\int_{c}^{d}f(a,y)dy \\
&&\left. -\left( \frac{r_{2}}{r_{2}+1}\right) \frac{1}{b-a}%
\int_{a}^{b}f(x,d)dx-\left( \frac{1}{r_{2}+1}\right) \frac{1}{b-a}%
\int_{a}^{b}f(x,c)dx\right\vert \\
&\leq &\frac{(b-a)(d-c)}{\left( r_{1}+1\right) \left( r_{2}+1\right) }\left( 
\frac{\left( 1+r_{1}^{2}\right) \left( 1+r_{2}^{2}\right) }{4\left(
r_{1}+1\right) \left( r_{2}+1\right) }\right) ^{1-\frac{1}{q}} \\
&&\times \left( \frac{MN\left\vert \frac{\partial ^{2}f}{\partial t\partial
\lambda }\right\vert ^{q}(a,c)+LN\left\vert \frac{\partial ^{2}f}{\partial
t\partial \lambda }\right\vert ^{q}(a,d)+KM\left\vert \frac{\partial ^{2}f}{%
\partial t\partial \lambda }\right\vert ^{q}(b,c)+KL\left\vert \frac{%
\partial ^{2}f}{\partial t\partial \lambda }\right\vert ^{q}(b,d)}{\left(
s+1\right) ^{2}\left( s+2\right) ^{2}}\right) ^{\frac{1}{q}}
\end{eqnarray*}%
for some fixed $r_{1},r_{2}\in \left[ 0,1\right] .$
\end{theorem}

\begin{proof}
From Lemma 3 and using the well-known Power-mean inequality, we can write%
\begin{eqnarray*}
&&\left\vert \frac{f\left( a,c\right) +r_{2}f\left( a,d\right) +r_{1}f\left(
b,c\right) +r_{1}r_{2}f\left( b,d\right) }{\left( r_{1}+1\right) \left(
r_{2}+1\right) }+\frac{1}{(b-a)(d-c)}\int_{a}^{b}\int_{c}^{d}f(x,y)dxdy%
\right. \\
&&-\left( \frac{r_{2}}{r_{2}+1}\right) \frac{1}{d-c}\int_{c}^{d}f(b,y)dy-%
\left( \frac{1}{r_{1}+1}\right) \frac{1}{d-c}\int_{c}^{d}f(a,y)dy \\
&&\left. -\left( \frac{r_{2}}{r_{2}+1}\right) \frac{1}{b-a}%
\int_{a}^{b}f(x,d)dx-\left( \frac{1}{r_{2}+1}\right) \frac{1}{b-a}%
\int_{a}^{b}f(x,c)dx\right\vert \\
&=&\frac{(b-a)(d-c)}{\left( r_{1}+1\right) \left( r_{2}+1\right) }\left(
\int_{0}^{1}\int_{0}^{1}\left\vert \left( \left( r_{1}+1\right) t-1\right)
\left( \left( r_{2}+1\right) \lambda -1\right) \right\vert dtd\lambda
\right) ^{1-\frac{1}{q}} \\
&&\times \left( \int_{0}^{1}\int_{0}^{1}\left\vert \left( \left(
r_{1}+1\right) t-1\right) \left( \left( r_{2}+1\right) \lambda -1\right)
\right\vert \left\vert \frac{\partial ^{2}f}{\partial t\partial \lambda }%
\left( tb+\left( 1-t\right) a,\lambda d+\left( 1-\lambda \right) c\right)
\right\vert ^{q}dtd\lambda \right) ^{\frac{1}{q}}.
\end{eqnarray*}%
Since $\left\vert \frac{\partial ^{2}f}{\partial t\partial \lambda }%
\right\vert ^{q}$ is $s-$convex function on the co-ordinates on $\Delta ,$
we can write for $t,\lambda \in \left[ 0,1\right] $%
\begin{eqnarray*}
&&\left\vert \frac{\partial ^{2}f}{\partial t\partial \lambda }\left(
tb+\left( 1-t\right) a,\lambda d+\left( 1-\lambda \right) c\right)
\right\vert ^{q} \\
&\leq &t^{s}\left\vert \frac{\partial ^{2}f}{\partial t\partial \lambda }%
\left( b,\lambda d+\left( 1-\lambda \right) c\right) \right\vert ^{q}+\left(
1-t\right) ^{s}\left\vert \frac{\partial ^{2}f}{\partial t\partial \lambda }%
\left( a,\lambda d+\left( 1-\lambda \right) c\right) \right\vert ^{q}
\end{eqnarray*}%
and%
\begin{eqnarray*}
&&\left\vert \frac{\partial ^{2}f}{\partial t\partial \lambda }\left(
tb+\left( 1-t\right) a,\lambda d+\left( 1-\lambda \right) c\right)
\right\vert ^{q} \\
&\leq &t^{s}\lambda ^{s}\left\vert \frac{\partial ^{2}f}{\partial t\partial
\lambda }\right\vert ^{q}(b,d)+t^{s}\left( 1-\lambda \right) ^{s}\left\vert 
\frac{\partial ^{2}f}{\partial t\partial \lambda }\right\vert ^{q}(b,c) \\
&&+\lambda ^{s}\left( 1-t\right) ^{s}\left\vert \frac{\partial ^{2}f}{%
\partial t\partial \lambda }\right\vert ^{q}(a,d)+\left( 1-\lambda \right)
^{s}\left( 1-t\right) ^{s}\left\vert \frac{\partial ^{2}f}{\partial
t\partial \lambda }\right\vert ^{q}(a,c)
\end{eqnarray*}%
hence, it follows that%
\begin{eqnarray}
&&\left\vert \frac{f\left( a,c\right) +r_{2}f\left( a,d\right) +r_{1}f\left(
b,c\right) +r_{1}r_{2}f\left( b,d\right) }{\left( r_{1}+1\right) \left(
r_{2}+1\right) }\right.  \notag \\
&&-\left( \frac{r_{2}}{r_{2}+1}\right) \frac{1}{d-c}\int_{c}^{d}f(b,y)dy-%
\left( \frac{1}{r_{1}+1}\right) \frac{1}{d-c}\int_{c}^{d}f(a,y)dy  \notag \\
&&-\left( \frac{r_{2}}{r_{2}+1}\right) \frac{1}{b-a}\int_{a}^{b}f(x,d)dx-%
\left( \frac{1}{r_{2}+1}\right) \frac{1}{b-a}\int_{a}^{b}f(x,c)dx  \label{6}
\\
&&\left. +\frac{1}{(b-a)(d-c)}\int_{a}^{b}\int_{c}^{d}f(x,y)dxdy\right\vert 
\notag \\
&\leq &\frac{(b-a)(d-c)}{\left( r_{1}+1\right) \left( r_{2}+1\right) }\left( 
\frac{\left( 1+r_{1}^{2}\right) \left( 1+r_{2}^{2}\right) }{4\left(
r_{1}+1\right) \left( r_{2}+1\right) }\right) ^{1-\frac{1}{q}}  \notag \\
&&\times \left( \int_{0}^{1}\int_{0}^{1}\left\vert \left( \left(
r_{1}+1\right) t-1\right) \left( \left( r_{2}+1\right) \lambda -1\right)
\right\vert \left\{ t^{s}\lambda ^{s}\left\vert \frac{\partial ^{2}f}{%
\partial t\partial \lambda }\right\vert ^{q}(b,d)\right. \right.  \notag \\
&&+t^{s}\left( 1-\lambda \right) ^{s}\left\vert \frac{\partial ^{2}f}{%
\partial t\partial \lambda }\right\vert ^{q}(b,c)+\lambda ^{s}\left(
1-t\right) ^{s}\left\vert \frac{\partial ^{2}f}{\partial t\partial \lambda }%
\right\vert ^{q}(a,d)  \notag \\
&&\left. \left. +\left( 1-\lambda \right) ^{s}\left( 1-t\right)
^{s}\left\vert \frac{\partial ^{2}f}{\partial t\partial \lambda }\right\vert
^{q}(a,c)\right\} dtd\lambda \right) ^{\frac{1}{q}}  \notag
\end{eqnarray}%
By a simple computation, one can see that%
\begin{eqnarray*}
&&\left( \int_{0}^{1}\int_{0}^{1}\left\vert \left( \left( r_{1}+1\right)
t-1\right) \left( \left( r_{2}+1\right) \lambda -1\right) \right\vert
\left\{ t^{s}\lambda ^{s}\left\vert \frac{\partial ^{2}f}{\partial t\partial
\lambda }\right\vert ^{q}(b,d)\right. \right. \\
&&+t^{s}\left( 1-\lambda \right) ^{s}\left\vert \frac{\partial ^{2}f}{%
\partial t\partial \lambda }\right\vert ^{q}(b,c)+\lambda ^{s}\left(
1-t\right) ^{s}\left\vert \frac{\partial ^{2}f}{\partial t\partial \lambda }%
\right\vert ^{q}(a,d) \\
&&\left. \left. +\left( 1-\lambda \right) ^{s}\left( 1-t\right)
^{s}\left\vert \frac{\partial ^{2}f}{\partial t\partial \lambda }\right\vert
^{q}(a,c)\right\} dtd\lambda \right) ^{\frac{1}{q}} \\
&=&\left( \frac{MN\left\vert \frac{\partial ^{2}f}{\partial t\partial
\lambda }\right\vert ^{q}(a,c)+LN\left\vert \frac{\partial ^{2}f}{\partial
t\partial \lambda }\right\vert ^{q}(a,d)+KM\left\vert \frac{\partial ^{2}f}{%
\partial t\partial \lambda }\right\vert ^{q}(b,c)+KL\left\vert \frac{%
\partial ^{2}f}{\partial t\partial \lambda }\right\vert ^{q}(b,d)}{\left(
s+1\right) ^{2}\left( s+2\right) ^{2}}\right) ^{\frac{1}{q}}
\end{eqnarray*}%
where $K,$ $L,$ $M$ and $N$ as in Theorem 11. By substituting these in (\ref%
{6}) and simplifying we obtain the required result.
\end{proof}

\begin{corollary}
(1) Under the assumptions of Theorem 13, if we choose $r_{1}=r_{2}=1,$ we
have%
\begin{eqnarray*}
&&\left\vert \frac{f\left( a,c\right) +f\left( a,d\right) +f\left(
b,c\right) +f\left( b,d\right) }{4}\right. \\
&&-\frac{1}{2}\left[ \frac{1}{d-c}\int_{c}^{d}\left[ f(b,y)+f(a,y)\right] dy+%
\frac{1}{b-a}\int_{a}^{b}\left[ f(x,d)+f(x,c)\right] dx\right] \\
&&\left. +\frac{1}{(b-a)(d-c)}\int_{a}^{b}\int_{c}^{d}f(x,y)dxdy\right\vert
\\
&\leq &\frac{(b-a)(d-c)}{4}\left( \frac{1}{4}\right) ^{1-\frac{1}{q}} \\
&&\times \left( \frac{MN\left\vert \frac{\partial ^{2}f}{\partial t\partial
\lambda }\right\vert ^{q}(a,c)+LN\left\vert \frac{\partial ^{2}f}{\partial
t\partial \lambda }\right\vert ^{q}(a,d)+KM\left\vert \frac{\partial ^{2}f}{%
\partial t\partial \lambda }\right\vert ^{q}(b,c)+KL\left\vert \frac{%
\partial ^{2}f}{\partial t\partial \lambda }\right\vert ^{q}(b,d)}{\left(
s+1\right) ^{2}\left( s+2\right) ^{2}}\right) ^{\frac{1}{q}}
\end{eqnarray*}%
(2) Under the assumptions of Theorem 13, if we choose $r_{1}=r_{2}=0,$ we
have%
\begin{eqnarray*}
&&\left\vert f\left( a,c\right) +\frac{1}{(b-a)(d-c)}\int_{a}^{b}%
\int_{c}^{d}f(x,y)dxdy\right. \\
&&\left. -\frac{1}{d-c}\int_{c}^{d}f(a,y)dy-\frac{1}{b-a}%
\int_{a}^{b}f(x,c)dx\right\vert \\
&\leq &(b-a)(d-c)\left( \frac{1}{4}\right) ^{1-\frac{1}{q}} \\
&&\times \left( \frac{MN\left\vert \frac{\partial ^{2}f}{\partial t\partial
\lambda }\right\vert ^{q}(a,c)+LN\left\vert \frac{\partial ^{2}f}{\partial
t\partial \lambda }\right\vert ^{q}(a,d)+KM\left\vert \frac{\partial ^{2}f}{%
\partial t\partial \lambda }\right\vert ^{q}(b,c)+KL\left\vert \frac{%
\partial ^{2}f}{\partial t\partial \lambda }\right\vert ^{q}(b,d)}{\left(
s+1\right) ^{2}\left( s+2\right) ^{2}}\right) ^{\frac{1}{q}}
\end{eqnarray*}
\end{corollary}

\begin{remark}
Under the assumptions of Theorem 13, if we choose $r_{1}=r_{2}=1$ and $s=1,$
we have an improvement for the inequality (\ref{1.5}).
\end{remark}

\end{document}